\documentclass[12pt]{article}
\title{Self-decomposability of Weak Variance Generalised Gamma Convolutions \footnote{This
research was partially supported by ARC grant DP160104737.}}

\author{Boris Buchmann\thanks{Research School of Finance, Actuarial Studies \& Statistics,
Australian National University,
ACT 0200,
Australia,
Phone (61-2) 6125 7296,
Fax (61-2) 6125 0087,
email: boris.buchmann@anu.edu.au
}
\and
Kevin W.~Lu\thanks{Research School of Finance, Actuarial Studies \& Statistics,
    Australian National University,
    ACT 0200,
    Australia,
    email: kevin.lu@anu.edu.au}
\and
Dilip B.~Madan\thanks{Robert H. Smith School of Business, University of Maryland, College Park, MD. 20742,
Maryland, USA, email: dbm@rhsmith.umd.edu}
}

\usepackage{amsmath,latexsym,amssymb,amsthm}
\usepackage{graphicx}	
\usepackage{hyperref}
\usepackage[all]{xy}

\newcommand{\rmd}{{\rm d}}
\newcommand{\rmi}{{\rm i}}

\newcommand{\eqd}{\stackrel{D}{=}}
\newcommand{\halmos}{\quad\hfill\mbox{$\Box$}}

\newcommand{\BB}{\mathbb{B}}
\newcommand{\CC}{\mathbb{C}}
\newcommand{\DD}{\mathbb{D}}
\newcommand{\EE}{\mathbb{E}}

\newcommand{\GG}{\mathbb{G}}

\newcommand{\RR}{\mathbb{R}}
\newcommand{\NN}{\mathbb{N}}

\newcommand{\PP}{\mathbb{P}}

\newcommand{\VV}{\mathbb{V}}
\newcommand{\mySS}{\mathbb{S}}

\newcommand{\myFF}{\mathbb{F}}

\newcommand{\AAAA}{\mathfrak{A}}
\newcommand{\DDDD}{\mathfrak{D}}
\newcommand{\EEEE}{\mathfrak{E}}
\newcommand{\HHHH}{\mathfrak{H}}

\newcommand{\KKKK}{\mathfrak{K}}

\newcommand{\hhhh}{\mathfrak{h}}

\newcommand{\PPP}{{\cal P}}

\newcommand{\GGG}{{\cal G}}

\newcommand{\SSS}{{\cal S}}
\newcommand{\TTT}{{\cal T}}
\newcommand{\UUU}{{\cal U}}
\newcommand{\VVV}{{\cal V}}

\newcommand{\XXX}{{\cal X}}

\newcommand{\GGC}{GGC}
\newcommand{\VGGC}{VGG}

\newcommand{\skal}[2]{\left\langle #1,#2\right\rangle}

\newcommand{\eins}{{\bf 1}}

\newcommand{\bfnull}{{\bf 0}}

\newcommand{\bfe}{{\bf e}}

\newcommand{\bfs}{{\bf s}}
\newcommand{\bft}{{\bf t}}
\newcommand{\bfx}{{\bf x}}
\newcommand{\bfy}{{\bf y}}

\newcommand{\bfz}{{\bf z}}
\newcommand{\bfd}{{\bf d}}
\newcommand{\bfu}{{\bf u}}
\newcommand{\bfv}{{\bf v}}
\newcommand{\bfB}{{\bf B}}
\newcommand{\bfT}{{\bf T}}
\newcommand{\bfX}{{\bf X}}
\newcommand{\bfY}{{\bf Y}}
\newcommand{\bfZ}{{\bf Z}}
\newcommand{\bfalpha}{\boldsymbol{\alpha}}
\newcommand{\bfdelta}{\boldsymbol{\delta}}

\newcommand{\bflambda}{\boldsymbol{\lambda}}
\newcommand{\bfmu}{\boldsymbol{\mu}}

\newcommand{\bfsigma}{\boldsymbol{\sigma}}

\newcommand{\bftheta}{\boldsymbol{\theta}}

\newcommand{\diag}{\mbox{diag}}

\newcommand{\wt}{\widetilde}
\newcommand{\tr}{\diamond}

\numberwithin{equation}{section}
\newtheorem{theorem}{Theorem}
\newtheorem{lemma}{Lemma}
\newtheorem{remark0}{\bf Remark}
\newenvironment{remark}{\begin{remark0}\em}{\end{remark0}\par}

\newtheorem{corollary}{Corollary}

\newtheorem{proposition}{Proposition}
\numberwithin{theorem}{section}
\numberwithin{proposition}{section}
\numberwithin{lemma}{section}
\numberwithin{corollary}{section}
\numberwithin{remark0}{section}
\numberwithin{definition}{section}

\begin{document}
\maketitle
\begin{abstract}
Weak variance generalised gamma convolution processes are multivariate Brownian motions weakly subordinated by
multivariate Thorin subordinators. Within this class, we extend a result from strong to weak subordination that
a driftless Brownian motion gives rise to a self-decomposable process. Under moment conditions on the underlying Thorin
measure, we show that this condition is also necessary. We apply our results
to some prominent processes such as the weak variance alpha-gamma process, and illustrate the necessity of our moment conditions in some cases.
\end{abstract}
\noindent
\begin{tabbing}
\end{tabbing}
\vspace{1cm}

\section{Introduction}\label{secintro}
An $n$-dimensional random vector $\bfX$ is called {\em self-decomposable}~($SD$) provided that for any $0\!<\!b\!<\!1$, there exists a random vector $\bfZ_b$, independent of $\bfX$, such that
\begin{align}\label{sddef1}
\bfX\;\eqd\;b\bfX+\bfZ_b\,.
\end{align}
Self-decomposable distributions occur as limits of sums of independent random vectors (see~\cite[Theorem~15.3]{s}) and as stationary distributions of multivariate L\'evy-driven Ornstein-Uhlenbeck processes~(see~\cite[Theorem~17.5]{s}). As such, self-decomposability has many important financial applications, such as modelling stochastic volatility~\cite{BNSh01}, stock returns \cite{BiKi02,Ma17}, and additional applications are reviewed in~\cite{Bi06,CGMY07}.
Self-decomposability was studied in the multivariate setting by Urbanik~\cite{Ur69}, who characterised their distributions in terms of a L\'evy-Khintchine representation, while Sato~\cite{sato80} derived a criterion often used to prove self-decomposability in terms of a representation of their L\'evy density in polar coordinates.

While it is clear that stable distributions, including the normal distribution and gamma distribution are $SD$, it was not until~\cite{Gr76} that self-decomposability was established for the Student's $t$-distribution. Thorin introduced generalised gamma convolutions ($GGC$), which are a class of self-decomposable distributions, to establish that both the Pareto distribution and the log-normal distribution are $SD$~\cite{Th77a,Th77b}.
His method was applied in~\cite{Ha79} to show that the generalised inverse Gaussian distribution is $GGC$, and thus $SD$ (see the survey article~\cite{JLY08} and the monographs~\cite{Bo92,SSV}).
Multivariate extensions of Thorin's results have been investigated in~\cite{BMS06,Bo09,BKMS16,BLMa,PS14}.

The $GGC^n$ class is the class of gamma distributions on rays of $[0,\infty)^n$ while being closed under convolution and distribution. They are infinitely divisible and the associated L\'evy processes are called Thorin subordinators. The question of whether the self-decomposability of the Thorin subordinator is inherited when subordinating Brownian motion is important and has been the subject of considerable research.

Grigelionis~\cite{Gr07} used {\em univariate} subordination of multivariate Brownian motion and {\em univariate} Thorin subordinators to construct the class of variance univariate generalised gamma convolutions ($VGG^{n,1}$), containing many prominent L\'evy processes used in mathematical finance~\cite{BiKi01,Eb01,MaSe90}. He showed that the corresponding $VGG^{n,1}$ distributions are $SD$ if $n\!=\!1$, or $n\!\geq\!2$ and the Brownian motion is driftless. However, if the Brownian motion is not driftless in the $n\geq 2$ case, then the distribution is not $SD$ under some moment conditions on the Thorin measure.% We rephrase his results in Proposition~\ref{propGrig} below.

In \cite{BKMS16}, the $VGG^{n,1}$ class was extended to the $VGG^{n,n}$ class of variance multivariate gamma convolutions using independent-component Brownian motion subordinated with {\em multivariate} Thorin subordinators. The $VGG^{n,n}$ class contains a number of recently introduced parametric classes of L\'evy processes~\cite{BKMS16,Gu13,LS10,Se08}. %\cite{Gu13,Se08} (see~\cite{BKMS16} and \cite{LS10}, their Subsection~2.6 and Section 3, respectively).
Barndorff-Nielsen, Pedersen and Sato \cite{BPS01} obtained sufficient conditions for the self-decomposability of processes formed by multivariate subordination. Applied to $VGG^{n,n}$ processes, this reduces to the Brownian motion subordinate being strictly stable, which is equivalent to driftlessness.

In~\cite{BLMa}, we introduced the overarching class of weak variance generalised gamma convolutions $WVGG^n$, using the notion of {\em weak} subordination, without the restriction that the Brownian motion has independent components. Weak subordination reduces to strong subordination if the subordinate has independent components or the subordinator has indistinguishable components, so $WVGG^n\!\supseteq\!VGG^{n,1}\cup VGG^{n,n}$. Subclasses of $WVGG^n$ processes, such as weak variance alpha-gamma processes, have been used in instantaneous portfolio theory~\cite{Ma17}, to model multivariate stock returns~\cite{BLMb,Mi17,MiSz17}, and the data analysis in \cite{BLMb} supports the hypothesis that log returns are self-decomposable.

In the present paper, we are concerned with the self-decomposability of $WVGG^n$ processes. We provide sufficient conditions and necessary conditions for this. In particular, we show that if the Brownian motion subordinate is driftless, then the $WVGG^n$ process is $SD$, while the converse holds under moment conditions on the Thorin measure, which extend and refine the results of Grigelions~\cite{Gr07}. In addition, we construct an example of a self-decomposable $WVGG^n$ process where the Brownian motion subordinate has nonzero drift, which illustrates the sharpness of our non-self-decomposability result. We obtain a complete characterisation of the self-decomposability of the weak variance alpha-gamma process provided the covariance of the Brownian motion is invertible.

Relatedly, there are two prominent generalisations of self-decomposability which were introduced by Urbanik. These are operator self-decomposability \cite{Ur72a} and the $L$ classes of nested operator self-decomposable distributions \cite{Ur72b}, which were further studied in \cite{sato80,SaYa85}. Operator self-decomposability allows for the previously mentioned correspondence with L\'evy-driven Ornstein-Uhlenbeck processes to be generalised to matrix-valued coefficients. Sufficient
 conditions for inclusion in these classes were studied in the context of multivariate subordination in \cite{BPS01}.

The remaining parts of the paper are organised as follows. In Section~\ref{secVGGC}, we introduce our notation, revise important properties of generalised variance gamma convolutions and gives some preliminary results. Section~\ref{secmain} contains our main results. Specifically, the sufficient conditions are given in Theorem~\ref{thmVGGCseldec}, while necessary conditions are given in Theorem~\ref{thmVGGCnotseldec}, and some equivalent conditions to one of the key moment conditions for non-self-decomposability are given in Proposition~\ref{propUUUAequiv}. In Section~\ref{secappl}, we apply our self-decomposability conditions to various examples, including the weak variance alpha-gamma process, and in Proposition~\ref{propsdcex}, we construct an $SD$ process where the Brownian motion subordinate has nonzero drift. We conclude in Remark~\ref{remconclude}. Section~\ref{secproofs} contains technical proofs.
\section{Notation and Preliminary Results}\label{secVGGC}
\noindent{\bf Notation.}~If $x,y\!\in\![-\infty,\infty]$ are extended real numbers, the minimum and maximum is denoted by $x\!\wedge\!y\!=\!\min\{x,y\}$ and $x\!\vee\!y\!=\!\max\{x,y\}$, so that $x\!=\!x^+\!-\!x^-$ decomposes $x$ into its positive and negative parts $x^+\!:=\!x\!\vee\!0$ and $x^-\!:=\!(-x)^+$, respectively.

Denote the principal branch of the complex logarithm by $\ln\!:\!\CC\backslash(-\infty,0]\!\to\!\CC$ so that, for $\rho\!\in\!\RR$, $w\!\in\!\CC$, $\Re w\!>\!0$, $w^\rho\!=\!e^{\rho\ln w}$ while the Sommerfeld integral representation of the modified Bessel function $K_{\rho}$ of second kind holds for $\rho\!\ge\!0$ and $w\!\in\mbox{dom}_\KKKK:=\{z\!\in\!\CC\!:\!\Re\!z\!>\!0,\,\Re\!z^2\!>\!0\}$, giving~(see~\cite[Equation~(3.471)--9]{GrRy96} and \cite{wa})
\begin{equation}\label{besselsommerfeld}\KKKK_{\rho}(w):=w^{\rho}K_{\rho}(w)
\,=\,2^{\rho-1}\int_0^\infty t^{\rho-1}\,
\exp\{-t\!-\!(w^2/4t)\}\,\rmd t\,.
\end{equation}
Let $\RR^n$ be $n$-dimensional Euclidean space, whose elements are row vectors $\bfx\!=\!(x_1,\dots,x_n)$, with canonical basis $\{\bfe_k\!:\!1\!\le\!k\!\le\!n\}$, and set
$\bfe\!:=\!(1,\dots,1)\!\in\!\RR^n$. Let $\bfx',\Sigma'$ denote the transpose of a vector $\bfx\in\RR^n$ and a matrix $\Sigma\in\RR^{n\times n}$, respectively. Let $\skal\bfx\bfy\!=\!\bfx \bfy'$ denote the Euclidean product with Euclidean norm $\|\bfx\|^2\!=\!\skal \bfx\bfx\!=\!\bfx\bfx'$, and set
$\skal \bfx \bfy_\Sigma\!:=\!\bfx \Sigma\bfy'$ and
$\|\bfx\|^2_\Sigma\!:=\!\skal{\bfx}{\bfx}_\Sigma$ for $\bfx,\bfy\!\in\!\RR^n$ and $\Sigma\!\in\!\RR^{n\times n}$ with determinant $|\Sigma|$.

Let $\DD\!:=\!\{\bfx\in\RR^n\!:\!\|\bfx\|\!\le\!1\}$ and $\mySS\!:=\!\{\bfx\in\RR^n\!:\!\|\bfx\|\!=\!1\}$ be the Euclidean unit ball and sphere centred at the origin, respectively.
If $A\!\subseteq\!\RR^n$, set $A_*\!:=\!A\!\backslash\!\{{\bf 0}\},A_{**}\!:=\!A\!\cap\!(\RR_*)^n,A_+\!:=\!A\!\cap\![0,\infty)^n,A_{++}\!:=\!A\!\cap\!(0,\infty)^n$ and  $A_{\le}\!:=\!\{\bfx\!=\!(x_1,\dots,x_n)\!\in\!A\!:x_1\le{ \dots}\le x_n\}$, and,
for $\bfy\!\in\!\RR^n$, let $\eins_A(\bfy)=\bfdelta_\bfy(A)$ denote the indicator function and the Dirac measure, respectively.

If $\bfx\!=\!(x_1,\dots,x_n)\!\in\!\RR^n$, set $\prod\bfx\!:=\!\prod_{k=1}^n x_k\!\in\!\RR$ while $\wedge\bfx\!:=\!(x_k\wedge x_l)_{1\le k,l\le n}\in\RR^{n\times n}$ is the matrix of component-wise minima. Let $*$ denote the Hadamard product of matrices. If $\bfmu\!=\!(\mu_1,\dots \mu_n)\!\in\!\RR^n$ and $\Sigma\!=\!(\Sigma_{kl})\!\in\!\RR^{n\times n}$, introduce~$\bfx\tr\bfmu\in\RR^n$ and~$\bfx\tr\Sigma\!\in\!\RR^{n\times n}$ as $\bfx\tr\bfmu\!:=\!\bfx*\bfmu\!:=\!(x_1\mu_1,\dots,x_n\mu_n)\in\RR^n$ and $\bfx\tr\Sigma\!:=\!(\wedge\bfx)*\Sigma$. Let the symmetric matrix $\Sigma^s\!:=\!(\Sigma\!+\!\Sigma')/2$ denote the symmetrisation of $\Sigma\!\in\!\RR^{n\times n}$. \\[1mm]
\noindent{\bf Matrix analysis.}~Oppenheim's and Hadamard's inequalities~(see~\cite[Theorems~3.7.5 and~3.6.3]{BR97}) state that $|\Sigma|\prod_{k=1}^n\Theta_{kk}
\!\le\!|\Sigma*\Theta|\!\le\!\prod_{k=1}^n(\Sigma_{kk}\Theta_{kk})$ for covariance matrices $\Sigma\!=\!(\Sigma_{kl}),\Theta\!=\!(\Theta_{kl})\!\in\!\RR^{n\times n}$.
If $\Sigma$ is a covariance matrix and $\bft\!\in\![0,\infty)^n$, Oppenheim's inequality ensures that $\bft\tr\Sigma$ is a covariance matrix.
If $\Sigma$ is an invertible covariance matrix and $\bft\!\in\!(0,\infty)^n$, Oppenheim's and Hadamard's inequalities ensure that $\bft\tr\Sigma$ is an invertible covariance matrix and, in addition,
\begin{equation}
\label{ineqOppHadutrSigma}
0\,<\,|\Sigma|\,\le\,\inf_{\bfu\in (0,\infty)^n}\frac{|\bfu\tr\Sigma|}{\prod \bfu}\,\le\,\sup_{\bfu\in (0,\infty)^n}\frac{|\bfu\tr\Sigma|}{\prod \bfu}\,\le\,\prod_{k=1}^n\Sigma_{kk}\,<\,\infty\,.
\end{equation}
Thus, for $\bfx,\bfmu\!\in\!\RR^n$, $\bfy\!\in\!\RR^n_*$, $\bfu\!\in\!(0,\infty)^n$, we can introduce well-defined quantities
\begin{align}\nonumber
\AAAA(\bfx,\bfu):=\AAAA_{\bfmu,\Sigma}(\bfx,\bfu)&:=\big\{(2\|\bfu\|^2\!+\!\|\bfu\tr\bfmu\|^2_{(\bfu\tr\Sigma)^{-1}})\|\bfx\|^2_{(\bfu\tr \Sigma)^{-1}}\big\}^{1/2}\;,\\
\DDDD(\bfx,\bfu):=\DDDD_{\Sigma}(\bfx,\bfu)&:=\|\bfx\|^{n}_{(\bfu\tr \Sigma)^{-1}}|\bfu\tr\Sigma|^{1/2}\,,\nonumber\\ \EEEE(\bfx,\bfu):=\EEEE_{\bfmu,\Sigma}(\bfx,\bfu)&:=\skal{\bfx}{\bfu\tr\bfmu}_{(\bfu\tr\Sigma)^{-1}}\,.\label{defADE}
\end{align}
Here $\AAAA(\bfx,\bfu),\DDDD(\bfy,\bfu)\in(0,\infty)$ and $\EEEE(\bfx,\bfu)\in\RR$ indicate ``exponent'', ``argument'' and ``denominator'' in~\eqref{defHHHH} below, respectively.

Next, we introduce the subset of $\RR^n$ where {\em uniform strict positivity (USP)} of $\EEEE$ holds by
\begin{equation}
\label{defVV}
\VV^+\,:=\,\VV_{\bfmu,\Sigma}^+
\,:=\,\big\{\bfv\in\RR^n:
\inf_{\bfu\in(0,\infty)^n}
\EEEE_{\bfmu,\Sigma}(\bfv,\bfu)\!>\!0
\big\}\, .
\end{equation}
We summarise properties of relating to $\EEEE$, $\AAAA$ and $\DDDD$ in Theorem~\ref{thmpos} (see Subsections~\ref{subsectoolsunifpositivity} and~\ref{subsecproofofthmpos} for a proof).
\begin{theorem}\label{thmpos}
Assume $\bfmu,\bfx\!\in\!\RR^n$, $\bfy\!\in\!\RR_*^n$, $\bfz\!\in\!\RR_{**}^n$ and an invertible covariance matrix $\Sigma\!\in\!\RR^{n\times n}$. Then:\\[1mm]
(i)~$\inf_{\bfu \in \mySS_{++}}\|\bfy\|_{(\bfu\tr\Sigma)^{-1}}>0$;\\[1mm]
(ii)~$\inf_{\bfu \in(0,\infty)^n}\DDDD(\bfz,\bfu)>0$;\\[1mm]
(iii)~$\sup_{\bfu\in\mySS_{++}}\|\bfu\tr\bfmu\|_{(\bfu\tr\Sigma)^{-1}}<\infty$;\\[1mm]
(iv)~$\sup_{\bfu\in(0,\infty)^n} |\EEEE(\bfx,\bfu)|<\infty$;\\[1mm]
(v)~$\VV^+$ is an open convex cone of $\RR^n$;\\[1mm]
(vi)~$\bfmu\!\in\!\VV^+$ when~$\bfmu\!\neq\!{\bf 0}$.\end{theorem}
\begin{remark}\label{remconverseinfsup}
    Theorem~\ref{thmpos}(i)--(iv) raises the question of whether the assertions $\sup_{\bfu \in \mySS_{++}}\|\bfx\|_{(\bfu\tr\Sigma)^{-1}}\!<\!\infty$, $\inf_{\bfu \in(0,\infty)^n}\DDDD(\bfy,\bfu)\!>\!0$, $\sup_{\bfu \in(0,\infty)^n}\DDDD(\bfz,\bfu)\!<\!\infty$, $\inf_{\bfu\in\mySS_{++}}\|\bfu\tr\bfmu\|_{(\bfu\tr\Sigma)^{-1}}\!>\!0$ and $\inf_{\bfu\in(0,\infty)^n} |\EEEE(\bfx,\bfu)|\!>\!0$ hold for $n\!\geq\!2$, $\bfmu,\bfx\!\in\!\RR^n$, $\bfy\!\in\!\RR_*^n$, $\bfz\!\in\!\RR_{**}^n$ and invertible covariance matrices $\Sigma\!\in\!\RR^{n\times n}$. These are all false in general. We obtain a counterexample for each of these assertions by letting $\Sigma$ be the identity matrix, and taking $\bfx\!\in\!\RR^n_{*}$, $\bfy\!=\!\bfe_1$, $\bfz\!\in\!\RR^n_{**}$, $\bfmu\!\in\!(\RR^n_{**})^C$, and $\bfx\!\in\!\RR^n$ such that $\skal{\bfx}{\bfmu}\!=\!0$, respectively. This illustrates a delicate balance in these quantities. We return to these facts in Remark~\ref{remadequiv1} below, while the importance of USP and the difficulty in proving it is discussed in Remark~\ref{remthmpos} below. \halmos
\end{remark}
\noindent{\bf L\'evy processes.}~We refer the reader to the monographs~\cite{Ap09,b,s} for necessary material on L\'e\-vy pro\-cesses while
our notation is adopted from~\cite{BKMS16,BLMa}.\\[1mm]
 The law of a L\'evy process \[\bfX\,=\,(X_1,\dots,X_n)\,=\,(\bfX(t))_{t\ge 0}\,=\,(X_1(t),\dots,X_n(t))_{t\ge 0}\] is determined by its characteristic function $\EE\exp\{\rmi\skal\bftheta{\bfX(t)}\}\!=\!\exp\{t\allowbreak\Psi_\bfX(\bftheta)\}$, $t\!\ge\!0$, $\bftheta\!\in\!\RR^n$, with associated L\'evy exponent $\Psi_\bfX\!=\!\Psi$, where
\begin{equation}\label{0.1}
    \Psi(\bftheta)\,=\,
    \rmi \skal {\bfmu}\bftheta\!-\!\frac 12\;\|\bftheta\|^2_{\Sigma}
    +\int_{\RR_*^n}\left(e^{\rmi\skal\bftheta \bfx}\!-\!1\!-\!\rmi\skal\bftheta \bfx \eins_\DD(\bfx)\right)\,\XXX(\rmd \bfx)\,.
\end{equation}
Here~$\bfmu\!=\!(\mu,\dots,\mu_n)\!\in\!\RR^n$, $\Sigma\!=\!(\Sigma_{kl})\!\in\!\RR^{n\times n}$ is a covariance matrix and $\XXX$ is a L\'evy measure, that is a nonnegative Borel measure on $\RR^n_*$ such that
$\int_{\RR^n_*}(1\wedge\|\bfx\|^2)\,\XXX(\rmd \bfx)\!<\!\infty$. We write $\bfX\!\sim\!L^n(\bfmu,\Sigma,\XXX)$ to denote that $\bfX$ is an $n$-dimensional L\'evy process with canonical triplet
$(\bfmu,\Sigma,\XXX)$. For $n$-dimensional L\'evy processes $\bfX$ and $\bfY$, $\bfX\!\eqd\!\bfY$ indicates that $\bfX$ and $\bfY$ are identical in law, that is $\Psi_{\bfX}\!=\!\Psi_{\bfY}$.\\[1mm]
\noindent{\bf Self-decom\-posability.}~An $n$-dimensional L\'evy process $\bfX$ is called {\em self-decom\-posable}~($SD^n$)
if the random vector $\bfX(1)$ is self-decomposable by \eqref{sddef1}. Gamma and Thorin subordinators below are self-decomposable~(see \cite[Equation (2.14)]{BMS06} and \cite{PS14}).

The self-decom\-posability of $\bfX$ is equivalent to its L\'evy measure having polar decomposition, for all Borel sets $A\!\subseteq\!\RR^n_*$,
\begin{equation}\label{eqselfdecompviapolar}
\XXX(A)=\int_{\mySS}\int_{(0,\infty)} \eins_{A}(r\bfs)k(\bfs,r) \frac{\rmd r}r\,\SSS(\rmd \bfs)\,,\end{equation}
where $\SSS$ is a finite Borel measure on the Euclidean sphere $\mySS$, $r\!\mapsto\!k(\bfs,r)$ is nonincreasing and nonnegative for all $\bfs\!\in\!\mySS$, and  $\bfs\!\mapsto\!k(\bfs,r)$ is Borel measurable for all $r\!>\!0$.\\[1mm]
\noindent{\bf Finite variation processes and subordinators.}~Sample paths of a L\'evy process $\bfX$ are of locally finite (or bounded) variation if and only if $\Sigma\!=\!0$ and
$\int_{\DD_*}\|\bfx\|\,\XXX(\rmd \bfx)\!<\!\infty$. In this case, $\bfd\!:=\!\bfmu\!-\!\int_{\DD_*} \bfx\,\XXX(\rmd \bfx)\!\in\!\RR^n$ denotes the drift of $\bfX$.
Particularly,~$\bfT\!=\!(T_1,\dots,T_n)\!\sim\!S^n(\bfd,\TTT)$ refers to an $n$-dimen\-sional subordinator, that is a L\'evy process with nondecreasing components, drift~$\bfd\!\in\![0,\infty)^n$ and L\'evy measure
$\TTT(([0,\infty)^n_*)^C)\!=\!0$.\\[1mm]
\noindent{\bf Gamma subordinator.}~Let $a,b\!>\!0$. A subordinator $G\!\sim\!\Gamma_S(a,b)$ is a {\em gamma subordinator} with shape $a$ and rate $b$ if and only if its marginal $G(t)\!\sim\!\Gamma(at,b)$, $t\!\geq\!0$,
is gamma distributed with shape $at$ and rate $b$. We have $G\!\sim\!S^1(0,\GGG_{a,b})$ with L\'evy measure $\GGG_{a,b}(\rmd g)\!:=\!\eins_{(0,\infty)}(g)a e^{-b g}\rmd g/g$ and Laplace exponent $-\ln\EE[\exp\{-\lambda G(t)\}]\!=\!at\ln\{(b\!+\!\lambda)/b\}$, $\lambda>-b$. If $a\!=\!b$, we refer to $G$ as a {\em standard} gamma subordinator, in short, $G\!\sim\!\Gamma_S(b)\!:=\!\Gamma_S(b,b)$. A gamma subordinator $G$ is a standard gamma subordinator if and only if $\EE[G(1)]\!=\!1$.\\[1mm]
\noindent{\bf Thorin subordinator.}~Our exposition follows~\cite{BKMS16}. A nonnegative Borel mea\-sure $\UUU$ on $[0,\infty)^n_*$ is called an $n$-dimensional {\em Thorin measure} provided
\begin{equation}\label{thorinmeasure}
    \int_{[0,\infty)_*^n}\;\left(1\!+\!\ln^-\|\bfu\|\right)\wedge\left(1\big/\|\bfu\|\right)\;\UUU(\rmd \bfu)\;<\;\infty\,.
\end{equation}
If $\bfd\!\in\![0,\infty)^n$ and $\UUU$ is a Thorin measure, we call an $n$-dimensional subordinator $\bfT$ a {\em Thorin subordinator}, in short $\bfT\!\sim\!\GGC_S^n(\bfd,\UUU)$, whenever, for all $t\!\ge\!0$, $\bflambda\!\in\![0,\infty)^n$, it has Laplace exponent
\[
-\ln \EE\exp\{-\skal{\bflambda}{\bfT(t)}\}\,=\,t\skal \bfd\bflambda+t\int_{[0,\infty)^n_*} \!\ln\big\{(\|\bfu\|^2+\skal{\bflambda}{\bfu})\big/\|\bfu\|^2\big\}\,\UUU(\rmd \bfu)\,.
\]
The distribution of a Thorin subordinator is uniquely determined by~$\bfd$ and~$\UUU$.

A driftless Thorin subordinator $T\!\sim\!\GGC^1_S(0,\UUU)$ is a gamma subordinator $T\!\sim\!\Gamma_S(a,b)$ if and only if $\UUU\!=\!a\bfdelta_b$. Alternatively, we may characterise the multivariate class of Thorin subordinators as the class of subordinators $G\bfalpha$ closed under convolutions and convergence in distribution, where $G$ is a gamma subordinator and $\bfalpha\!\in\![0,\infty)^n$~\cite{PS14}.\\[1mm]
\noindent{\bf Brownian motion.}~Throughout, we let $\bfB\!=\!(B_1,\dots,B_n)\!\sim\!BM^n(\bfmu,\Sigma)\!:=\!L^n(\bfmu,\Sigma,0)$ denote $n$-dimensional Brownian motion $\bfB$ with linear drift $\EE[\bfB(t)]\!\allowbreak=\!t\bfmu$ and covariance matrix Cov$(\bfB(t))\!=\!t\Sigma$, $t\!\ge\!0$.\\[1mm]
{\bf Strongly subordinated Brownian motion.}~Let $\bfT\!=\!(T_1,\dots,T_n)$ be an $n$-dimensional subordinator, and $\bfB\!=\!(B_1,\dots,B_n)$ be an $n$-dimensional Brownian motion. The {\em strong subordination of $\bfB$ by $\bfT$} is the $n$-dimensional process $\bfB\!\circ\!\bfT$ defined by
\begin{align}\label{strongsubord}
    (\bfB\circ\bfT)(t):=(B_1(T_1(t)),\dots,B_n(T_n(t)))\,,\quad t\ge 0\,.
\end{align}
For independent $\bfB$ and $\bfT$, $\bfB\!\circ\!\bfT$ is a L\'evy process if $\bfT\!\equiv\!T_1\bfe$ has indistinguishable components~\cite{BS10,R065,s,Zo58} or
$\bfB$ has independent components $B_1,,\dots,B_n$~\cite{BPS01}, otherwise it may not be~(see~\cite[Proposition~3.9]{BLMa}).\\[1mm]
\noindent{\bf Weakly subordinated Brownian motion.}~If $\bfB\!\sim\!BM^n(\bfmu,\Sigma)$ is Brownian motion and $\bfT\!\sim\!S^n(\bfd,\TTT)$ is a subordinator, we refer to an $n$-dimensional L\'evy process $\bfX\!\eqd\!\bfB\!\odot\!\bfT$ as {\em $\bfB$ weakly subor\-dinated by $\bfT$} provided
\[\bfX\,\sim\,L^n\big(\bfd\tr\bfmu\!+\!\int_{\DD_*}\!\bfx\XXX(\rmd\bfx),\,\bfd\tr\Sigma,\XXX\big)\,,\]
where $\XXX(\rmd\bfx)\!=\!\int_{[0,\infty)_*^n}\PP(\bfB(\bft)\!\in\!\rmd\bfx)\TTT(\rmd\bft)$ is a L\'evy measure~(see~\cite[Proposition~3.2]{BLMa}).
Alternatively, $\bfX\!\eqd\!\bfB\!\odot\!\bfT$ if and only if its characteristic exponent is (see~\cite[Proposition~3.1]{BLMa}), $\bftheta\!\in\!\RR^n$,
\[\Psi_\bfX(\bftheta)=\rmi\skal\bftheta{\bfd\tr\bfmu}-\frac 12\|\bftheta\|^2_{\bfd\tr\Sigma}+\int_{[0,\infty)^n_*}\!
\big(\exp\big\{\rmi\skal\bftheta{\bft\tr\bfmu}-\frac 12\|\bftheta\|^2_{\bft\tr\Sigma}\big\}-1\big)\TTT(\rmd\bft).\]
If $\bfB$ and $\bfT$ are independent, $\bfB\!\odot\!\bfT\!\eqd\!\bfB\!\circ\!\bfT$ provided $\bfT$ has indistinguishable components or
$\bfB$ has independent components~(see \cite[Proposition~3.3]{BLMa}). In contrast to strong subordination, using weak subordination, we stay within the framework of L\'evy processes while allowing $\bfB$
to have general covariance matrices $\Sigma$. For detailed accounts we refer to~\cite{BLMa,Lu18}.\\[1mm]
\noindent{\bf Variance gamma.}~For a Brownian motion $\bfB\!\sim\!BM^n(\bfmu,\Sigma)$ independent of a gamma subordinator $G\!\sim\!\Gamma_S(b)$, we call $\bfX\!\sim\!VG^n(b,\bfmu,\Sigma)$
a {\em variance gamma process}~\cite{MaSe90} provided $\bfX\!\eqd\!\bfB\!\circ\!(G\bfe)$.\\[1mm]
\noindent{\bf Strong variance univariate} $\boldsymbol{GGC}${\bf .}~Let $d\!\ge\!0$, $\bfmu\!\in\!\RR^n$, $\Sigma\!\in\!\RR^{n\times n}$ be an {\em arbitrary} covariance matrix and $\UUU$ be a {\em univariate} Thorin measure.  We call a L\'evy process $\bfX\!\sim\! \VGGC^{n,1}(d,\bfmu,\Sigma,\UUU)$ an $n$-dimen\-sional {\em strong variance univariate generalised gamma convolution process}~\cite{BKMS16,Gr07} provided $\bfX\!\eqd\!\bfB\!\circ\!(T\bfe)$, where $\bfB\!\sim\!BM^n(\bfmu,\Sigma)$ and $T\!\sim\!\GGC^1_S(d,\UUU)$ are independent.\\[1mm]
\noindent\noindent{\bf Strong variance multivariate} $\boldsymbol{GGC}${\bf .}~Let $\bfd\!\in\![0,\infty)^n$, $\bfmu\!\in\!\RR^n$, $\Sigma\!\in\!\RR^{n\times n}$ be a {\em diagonal} covariance matrix and $\UUU$ be an $n$-dimensional Thorin measure. We call a L\'evy process $\bfX\!\sim\!\VGGC^{n,n}(\bfd,\bfmu,\Sigma,\UUU)$ an $n$-dimensional {\em strong variance multivariate generalised gamma convolution process}~\cite{BKMS16} provided $\bfX\!\eqd\!\bfB\!\circ\!\bfT$, where $\bfB\!\sim\!BM^n(\bfmu,\Sigma)$ and $\bfT\!\sim\!\GGC^n_S(\bfd,\UUU)$ are independent.\\[1mm]
\noindent{\bf Weak variance generalised gamma convolutions.}~Let $\bfd\!\in\![0,\infty)^n$, $\bfmu\!\in\!\RR^n$, $\Sigma\!\in\!\RR^{n\times n}$ be an {\em arbitrary} covariance matrix and $\UUU$ be an $n$-dimensional Thorin measure. A L\'evy process $\bfX\!\sim\! W\VGGC^{n}(\bfd,\bfmu,\Sigma,\UUU)$ is called an $n$-dimensional {\em weak variance generalised gamma convolution process}~\cite{BLMa} provided $\bfX\!\eqd\!\bfB\!\odot\!\bfT$, where $\bfB\!\sim\!BM^n(\bfmu,\Sigma)$ and $\bfT\!\sim\!\GGC^n_S(\bfd,\UUU)$.

The notation $\VGGC^{n}\!=\!W\VGGC^{n}$ has also been used~\cite{BLMa,Lu18}.
If $\bfX\!\sim\!W\VGGC^{n}(\bfd,\bfmu,\Sigma,\UUU)$, the associated characteristic exponent is~(see~\cite[Theorem~4.1]{BLMa}), $\bftheta\in\RR^n$,
\begin{eqnarray}
    \Psi_\bfX(\bftheta)&=&
    \rmi \skal{\bfd\tr\bfmu}{\bftheta}-\frac 12 \|\bftheta\|^2_{\bfd\tr\Sigma}\label{GVGcharexpo}
    \\
    &&{}-\int_{[0,\infty)^n_*}\!\ln\big\{(\|\bfu\|^2
    -\rmi\skal{\bfu\tr\bfmu}{\bftheta}+\frac 12\|\bftheta\|^2_{\bfu\tr\Sigma})\big/\|\bfu\|^2\big\}\UUU(\rmd\bfu)\,.\nonumber
\end{eqnarray}
Note $\bfX\!\sim\!W\VGGC^{n}(\bfd,\bfmu,\Sigma,\UUU)$ is a
$\VGGC^{n,1}$ process if and only if $\bfd\!=\!d\bfe$ and $\UUU\!=\!\int_{0+}^\infty\bfdelta_{v\bfe/\|\bfe\|^2}\,\UUU_0(\rmd v)$
for some $d\!\ge\!0$ and a univariate Thorin measure $\UUU_0$ (see~\cite[Remark~4.2]{BLMa}). If $n\!=\!1$, our notation collapses into one inclusion $VG^1\!\subseteq\!\VGGC^{1,1}\!=\!W\VGGC^1$. If $n\!\ge\!2$, the $\VGGC^{n,1}$ class and $\VGGC^{n,n}$ class appear complementary~\cite{BKMS16}; if $n\!\ge\!1$, we have the chain of inclusions $VG^n\!\subseteq\!\VGGC^{n,1}\!\cup\!\VGGC^{n,n}\!\subseteq\!W\VGGC^n$~\cite{BLMa}.

Recall~\eqref{defADE}. If $\bfy\!\in\!\RR_{*}^n$, setting $c_n\!:=\!2/(2\pi)^{n/2}$, introduce
the function $w\!\mapsto\!\HHHH_\bfy(w)$, $w\!\in\!\mbox{dom}_\KKKK$, by
\begin{eqnarray}\label{defHHHH}
\HHHH_\bfy(w)\,:=\,c_n\int_{(0,\infty)^n}\exp\{w\EEEE(\bfy,\bfu)\}\,
\KKKK_{n/2}\{w \AAAA(\bfy,\bfu)\}\,\frac{\UUU(\rmd\bfu)}{\DDDD(\bfy,\bfu)}\,.
\end{eqnarray}
In~\eqref{XXXviaHHHH} below, $\HHHH_\bfs(r)/r$, $r\!>\!0$, $\bfs\!\in\!\mySS_{**}$, occurs as the density of the $WVGG^n$ L\'evy measure when restricted to $\RR_{**}^n$ in Euclidean polar coordinates, and as such,
it is the multivariate analogue of the $\tilde k$ function in~\cite{Gr07} (also see~\eqref{eqselfdecompviapolar}).

Next, we summarise some of its properties in Theorem~\ref{thmHAEDformula} (see Subsection~\ref{subsecproofthmHAEDformula} for a proof).
\begin{theorem}\label{thmHAEDformula} Let $\bfX\!\sim\!W\VGGC^n(\bfd,\bfmu,\Sigma,\UUU)$ with L\'evy measure $\XXX$, $n\!\ge\!2$, $|\Sigma|\neq 0$ and $\bfz\in\RR_{**}^n$. Then we have:\\[1mm]
(i)~$w\mapsto\HHHH_\bfz(w)$ is holomorphic on $\mbox{{\em dom}}_\KKKK$ while $\HHHH_\bfz((0,\infty))\subseteq[0,\infty)$, and its complex derivative is computed under the integral as
\begin{eqnarray}\label{partialrHHHH}
\lefteqn{\hspace{-5mm}\partial_w\HHHH_\bfz(w)\,=\,c_n\int_{(0,\infty)^n}\EEEE(\bfz,\bfu)\exp\{w\EEEE(\bfz,\bfu)\}\,
\KKKK_{n/2}\{w\AAAA(\bfz,\bfu)\}\,\frac{\UUU(\rmd\bfu)}{\DDDD(\bfz,\bfu)}}&&\\
&&-c_nw\,\int_{(0,\infty)^n}\AAAA^2(\bfz,\bfu)\exp\{w\EEEE(\bfz,\bfu)\}\,
\KKKK_{(n-2)/2}\{w\AAAA(\bfz,\bfu)\}\frac{\UUU(\rmd\bfu)}{\DDDD(\bfz,\bfu)}\,.\nonumber\end{eqnarray}
(ii)~If $\int_{(0,\infty)^n}\!\AAAA(\bfz,\bfu)\UUU(\rmd\bfu)/\DDDD(\bfz,\bfu)\!<\!\infty$, then $\int_{(0,\infty)^n}\!
|\EEEE(\bfz,\bfu)|\UUU(\rmd\bfu)/\DDDD(\bfz,\bfu)\!<\!\infty$ and
\begin{equation}\label{HHHHdiffat0}
\partial_r\HHHH_\bfz(0+)=c_n\,2^{(n-2)/2}\Gamma(n/2)\,\int_{(0,\infty)^n}\EEEE(\bfz,\bfu)\,\frac{\UUU(\rmd\bfu)}{\DDDD(\bfz,\bfu)}\,.\end{equation}
(iii)~If $A\subseteq\RR_{**}^n$ is a Borel set, we have
\begin{equation}\label{XXXviaHHHH}
\XXX(A)\,=\,\int_{\mySS_{**}}\int_{0+}^{\infty}\eins_{A}(r\bfs)\,\HHHH_\bfs(r)\,\frac{\rmd r}r\,\rmd\bfs\,,
\end{equation}
where $\rmd\bfs$ is the Lebesgue surface measure on $\mySS$.
\end{theorem}
\noindent\noindent{\bf Generalised variance gamma convolutions.}~The class of {\em generalised variance gamma convolutions} ($GVGC^n$) is the class of $VG^n(b,\bfmu,\Sigma)$ processes closed under convolution and convergence in distribution. Similarly, the class of $VG^n(b,\bfnull,\Sigma)$ processes closed under convolution and convergence in distribution is denoted $GVGC^n_{\bfnull}$. A process $\bfX\!\sim\!BM^n\!+\!GVGC_{\bfnull}^n$ if and only if $\bfX\!\eqd \!\bfB\!+\!\bfY$, where $\bfB\!\sim\!BM^n(\bfmu,\Sigma)$ and $\bfY\!\sim\!GVGC^n_{\bf0}$ are independent.\\[1mm]
{\bf Driftless subordination classes.}~For the process classes defined, we use the subscripts $\bfnull$ and $*$ to denote the restriction to the cases where $\bfmu=\bfnull$ and $\bfmu\neq\bfnull$, respectively.

\noindent{\bf Self-decom\-posability, revisited.}~The self-decomposability of the $\VGGC^{n,1}$ class was investigated in~\cite[Proposition~3]{Gr07}, and this result is stated here as Proposition~\ref{propGrig}. If $n\!=\!1$, the result was shown in~\cite[Section~3]{Ha79}.
\begin{proposition}\label{propGrig}Assume $\bfX\!\sim\!\VGGC^{n,1}(d,\bfmu,\Sigma,\UUU)$.\\
    (i)~If $n\!=\!1$, or $n\!\ge\! 2$ and $\bfmu\!=\!{\bf 0}$, then $\bfX\sim SD^n$.\\
    (ii)~If $n\!=\!2$, $|\Sigma|\!\neq\!0$, $\bfmu\!\neq\!{\bf 0}$ and $0<\int_{(0,\infty)} u^2\,\UUU(\rmd u)<\infty$, or $n\!\ge\!3$, $|\Sigma|\!\neq\!0$, $\bfmu\!\neq\!{\bf 0}$ and $0<\int_{(0,\infty)}u\,\UUU(\rmd u)<\infty$, then $\bfX\not\sim SD^n$.\end{proposition}
In context of the $\VGGC^{n,n}$ class, we can deduce the following sufficient condition for self-decomposability from~\cite[Theorem 6.1]{BPS01}.
\begin{proposition}\label{propBPS}Assume $\bfX\!\sim\!\VGGC^{n,n}_\bfnull$, then $\bfX\!\sim\!SD^n$.\end{proposition}
\section{Main Results}\label{secmain}
\noindent{\bf Sufficient conditions.}~As a generalisation of Propositions~\ref{propGrig}(i) and~\ref{propBPS}, we give a sufficient condition to ensure that a $W\VGGC^n$ process is self-decomposable~(see Subsection~\ref{subsecProofsthmselfdec} for a proof).
\begin{theorem}\label{thmVGGCseldec} Always, $VG^1\!\subseteq\!\VGGC^{1,1}\!=\!W\VGGC^1\!\subseteq\!SD^1$. Otherwise, if $n\!\ge\!2$, then the implications (i)$\Leftarrow$(ii)$\Leftarrow$(iii) hold,
    where:\\[1mm]
    (i)~$\bfX\!\sim\!SD^n$;\\[1mm]
    (ii)~$\bfX\!\sim\!BM^n\!+\!GVGC^n_{\bf0}$;\\[1mm]
    (iii)~$\bfX\!\sim\!W\VGGC^n_{\bf0}$.
\end{theorem}
\noindent{\bf Necessary conditions.}~As a generalisation of Proposition~\ref{propGrig}(ii), we give a sufficient condition to ensure that a $W\VGGC^n$ process is not self-decomposable (see Subsection~\ref{subsecProofthmsnotseldec} for a proof). Recall~\eqref{defHHHH}.
\begin{theorem}\label{thmVGGCnotseldec} If $\bfX\!\sim\!W\VGGC^n(\bfd,\bfmu,\Sigma,\UUU)$, $n\!\ge\!2$ and
$|\Sigma|\!\neq \!0$, then the implications (i)$\Leftarrow$(ii)$\Leftarrow$(iii)$\Leftarrow$(iv)$\Leftarrow$(v)$\Leftarrow$(vi)$\Leftarrow$(vii) hold,
where:\\[1mm]
(i)~$\bfX\!\not\sim\!SD^n$;\\[1mm]
(ii)~there exist a Borel set $\BB\!\subseteq\!\mySS_{**}$ of strictly positive $(n\!-\!1)$-dimensional Lebesgue surface mea\-sure such that, for all $\bfs\!\in\!\BB$,
$r\!\mapsto\!\HHHH_\bfs(r)$ is strictly increasing at some point $r_0\!\in\!(0,\infty)$;\\[1mm]
(iii)~there exist a Borel set $\BB\!\subseteq\!\mySS_{**}$ of strictly positive $(n\!-\!1)$-dimensional Lebesgue surface mea\-sure such that, for all $\bfs\!\in\!\BB$, the right-hand limit $\partial_r\HHHH_\bfs(0+)$ exists and is strictly positive;\\[1mm]
(iv)~there exist a Borel set $\BB\!\subseteq\!\mySS_{**}$ of strictly positive $(n\!-\!1)$-dimensional Le\-besgue surface mea\-sure such that, for all $\bfs\!\in\!\BB$,
\begin{equation}\label{UUUAintegrab}
\int_{(0,\infty)^n}\,\AAAA(\bfs,\bfu)\,\frac{\UUU(\rmd\bfu)}{\DDDD(\bfs,\bfu)}\;<\;\infty
\end{equation}
and
\begin{equation}\label{meanpositvity}
\int_{(0,\infty)^n}\EEEE(\bfs,\bfu)\,\frac{\UUU(\rmd\bfu)}{\DDDD(\bfs,\bfu)}\;>\,0\;;
\end{equation}
(v)~$\bfmu\!\neq\!{\bf 0}$ and
\begin{equation}\label{UUUAintegrabstronger}
0\;<\;\int_{(0,\infty)^n}(1\!+\!\|\bfu\|^{1/2})\,\left(\frac{\|\bfu\|^{n}\,}{\prod\bfu}\right)^{1/2}\UUU(\rmd\bfu)\;<\;\infty\,;\end{equation}
(vi)~$\bfmu\!\neq\!{\bf 0}$ and $\UUU((0,\infty)^n\cap\cdot)\!=\!\sum_{k=1}^m\int_{(0,\infty)} \bfdelta_{v\bfu_k}(\cdot)\,\UUU_k(\rmd v)$ for some $\bfu_k\!\in\!(0,\infty)^n$ and univariate Thorin measures $\UUU_k$, where
\begin{equation}\label{defraymoments}
0\,<\,\int_{(0,\infty)} v^{1/2}\,\UUU_k(\rmd v)\,<\,\infty\,,\qquad 1\le k\le m\,,\; m\ge 1\,.
\end{equation}
(vii)~$\bfmu\!\neq\!{\bf 0}$ and $\UUU((0,\infty)^n)>0$ for a finitely supported Thorin measure $\UUU$.
\end{theorem}
\begin{remark} While~\eqref{UUUAintegrab} and~\eqref{meanpositvity} show delicate dependencies on $\bfs$, $\bfmu$ and $\Sigma$, we can replace these with the more robust condition~\eqref{UUUAintegrabstronger}.
Clearly,~\eqref{meanpositvity} vanishes for $\bfmu\!=\!{\bf 0}$, which is consistent with Theorem~\ref{thmVGGCseldec}. The assumption~\eqref{UUUAintegrab}
ensures the existence of the finite right-hand limit of the derivative at the origin. If $\bfd\!=\!\bfnull$, then~\eqref{UUUAintegrabstronger} entails that $\bfX$ has paths of finite variation~(see~\cite[Proposition 2.6(a)]{BKMS16} and~\cite[Proposition~4.1]{BLMa}).~\halmos
\end{remark}
\begin{remark}\label{remthmpos}~%The importance of USP (see Theorem~\ref{thmpos}(v)--(vi)) is made clear by~\eqref{meanpositvity}.
    From Theorem~\ref{thmVGGCnotseldec}(iv), we see that USP (see Theorem~\ref{thmpos}(v)--(vi)) is important because in ensures that \eqref{meanpositvity} holds.

    Assume $n\!\ge\!1$, $\bfmu\!\in\!\RR^n_*$. We can easily verify USP in a limited number of cases. For Grigelionis' $VGG^{n,1}$ class~\cite{Gr07}, note $\EEEE_{\bfmu,\Sigma}(\bfmu,u\bfe)\!=\!\|\bfmu\|^2_{\Sigma^{-1}}\!>\!0$ for $u\!>\!0$, and for the $\VGGC^{n,n}$ class~\cite{BKMS16},
    note $\EEEE_{\bfmu,\Sigma}(\bfmu,\bfu)\!=\!\|\bfmu\|^2_{\Sigma^{-1}}\!>\!0$ for $\bfu\!\in\!(0,\infty)^{n}$ as $\Sigma$ is diagonal.

    But if $\Sigma$ is not a diagonal matrix, we have $\EEEE_{\bfmu,\Sigma}(\bfmu,u\bfalpha)\!=\!\skal\bfmu{\bfalpha\tr\bfmu}_{(\bfalpha\tr\Sigma)^{-1}}$ for $u\!>\!0$, $\bfalpha\!\in\!(0,\infty)^n$, but the strict positivity of $\skal\bfmu{\bfalpha\tr\bfmu}_{(\bfalpha\tr\Sigma)^{-1}}$ is not obvious (this case include the $WVAG^n$ process introduced below). This illustrates that how elusive USP is for overarching $W\VGGC^n$ class. \halmos
    %This illustrates that USP is elusive for the $WVAG^n$ class, let alone the overarching $W\VGGC^n$ class.\halmos
\end{remark}
For $VGG^{n,1}$ and $VGG^{n,n}$ processes, we obtain the following corollaries as implications of Theorems~\ref{thmVGGCseldec} and~\ref{thmVGGCnotseldec}. The result for $VGG^{n,1}$ processes is a refinement of Proposition~\ref{propGrig}.
\begin{corollary}\label{corVGGn1}Assume $\bfX\!\sim\!\VGGC^{n,1}(d,\bfmu,\Sigma,\UUU)$.\\
(i)~If $n\!=\!1$, or $n\!\ge\! 2$ and $\bfmu\!=\!{\bf 0}$, then $\bfX\!\sim\!SD^n$.\\
(ii)~If $n\!\ge\!2$, $|\Sigma|\!\neq\!0$,~$\bfmu\!\neq\!{\bf 0}$ and $0<\int_{(0,\infty)} u^{1/2}\,\UUU(\rmd u)\!<\!\infty$,~then $\bfX\!\not\sim\!SD^n$.
\end{corollary}
\begin{corollary}\label{corVGGnn} Assume~$\bfX\!\sim\!\VGGC^{n,n}(\bfd,\bfmu,\Sigma,\UUU)$ for a diagonal matrix $\Sigma$.\\
(i)~If $n\!=\!1$ or $n\!\ge\! 2$ and $\bfmu\!=\!{\bf 0}$, then $\bfX\!\sim\!SD^n$.\\
(ii)~If $n\!\ge\!2$, $|\Sigma|\!\neq\!0$, $\bfmu\!\neq\!{\bf 0}$ and~\eqref{UUUAintegrabstronger} holds, then $\bfX\!\not\sim\!SD^n$.
\end{corollary}

The following proposition gives equivalent conditions for \eqref{UUUAintegrab} (see Subsection~\ref{subsecProofUAequi} for a proof).

\begin{proposition}\label{propUUUAequiv}
    Let $n\!\ge\!2$, $\Sigma\!\in\!\RR^{n\times n}$, $|\Sigma|\!\neq \!0$ and $\bfs\!\in\!\mySS_{**}$.\\
    (i)~If $\UUU((a\DD_*)_{++})\!=\!0$ for some $a\!>\!0$, then \eqref{UUUAintegrab} is equivalent to
    \begin{equation}\label{UUUAequiv}
    %\int_{(0,\infty)^n} \frac{(1+\|\bfu\|)^{1/2}\|\bfu\|^{1/2}}{\|\bfs\|^{n-1}_{(\bfu\tr\Sigma)^{-1}}(\prod\bfu)^{1/2}}\,\UUU(\rmd\bfu)\,<\,\infty\,.
    \int_{((a\DD)^C)_{++}} \frac{\|\bfu\|}{\|\bfs\|^{n-1}_{(\bfu\tr\Sigma)^{-1}}(\prod\bfu)^{1/2}}\,\UUU(\rmd\bfu)\,<\,\infty\,.
    \end{equation}
    (ii)~If $\UUU((\mySS_{++})^C)\!=\!0$, then \eqref{UUUAintegrab} is equivalent to
    \begin{equation}\label{UUUAequivcirc}
    \int_{\mySS_{++}} {\|\bfs\|^{1-n}_{(\bfu\tr\Sigma)^{-1}}\big(\prod\bfu\big)^{-1/2}}\,\UUU(\rmd\bfu)\,<\,\infty\,.
    \end{equation}
    (iii)~If $\UUU((0,\infty)^n\cap\cdot)\!=\!\sum_{k=1}^m\int_{(0,\infty)} \bfdelta_{v\bfu_k}(\cdot)\,\UUU_k(\rmd v)$ for some $\bfu_k\!\in\!(0,\infty)^n$ and univariate Thorin measures $\UUU_k$, then \eqref{UUUAintegrab} is equivalent to
    \begin{equation}\label{defraymomentsequiv}
    \int_{(1,\infty)} v^{1/2}\,\UUU_k(\rmd v)\,<\,\infty\,,\qquad 1\le k\le m\,,\; m\ge 1\,.
    \end{equation}
\end{proposition}
\begin{remark}\label{remadequiv1}
    Fix any $\bfs\!\in\!\mySS_{**}$. Note~\eqref{UUUAintegrabstronger} and~\eqref{defraymoments} are sufficient for~\eqref{UUUAintegrab} to hold, and Theorem~\ref{thmpos}(ii) shows that $\int_{(0,\infty)^n}\,\AAAA(\bfs,\bfu)\,{\UUU(\rmd\bfu)}\!<\!\infty$ is also sufficient. Proposition~\ref{propUUUAequiv} give equivalent conditions to~\eqref{UUUAintegrab} under various assumptions on the Thorin measure, while Remark~\ref{remarkwhynogenadequiv} below explains why cannot find an equivalent condition for a general Thorin measure. Remark~\ref{remconverseinfsup} tells us that $\int_{(0,\infty)^n}\,\AAAA(\bfs,\bfu)\,{\UUU(\rmd\bfu)}\!<\!\infty$ may not be an equivalent condition as $\sup_{\bfu \in(0,\infty)^n}\DDDD(\bfs,\bfu)\!=\!\infty$ can occur, and that we cannot necessarily remove the terms $\|\bfs\|^{n-1}_{(\bfu\tr\Sigma)^{-1}}$ or $(\prod\bfu)^{1/2}$ from~\eqref{UUUAequiv}--\eqref{UUUAequivcirc} while remaining an equivalent condition as their infimums over $\bfu\!\in\!(0,\infty)^n$ are 0.
\end{remark}
\section{Examples}\label{secappl}
To give examples, we apply our sufficient and necessary conditions to various classes of $WVGG^n$ processes, which we review below. In our examples, the Thorin measure is often supported on a finite number of points~(see~\cite[Subsection~2.5]{BKMS16}), or more generally, on a finite union of rays~\cite{LS10}. The latter case includes multivariate Thorin subordinators that model common and idiosyncratic time changes such as
the alpha-gamma subordinators, which are often of interest in financial applications~\cite{LS10}.\\[1mm]
\noindent{\bf Alpha-gamma subordinator.}~Assume $n\!\ge\!2$. Let $a\!>\!0$, $\bfalpha\!=\!(\alpha_1,\dots,\alpha_n)\!\in\!(0,\infty)^n$
such that $a\alpha_k\!<\!1$ for $k=1,\dots,n$. Introduce $\beta_k\!:=\!(1\!-\!a\alpha_k)/\alpha_k$, and let $G_0,\dots,G_{n}$ be independent gamma subordinators such that $G_{0}\!\sim\!\Gamma_S(a,1)$, $G_k\!\sim\!\Gamma_S(\beta_k,1/\alpha_k)$, $1\!\le\!k\!\le\!n$.

\begin{figure}[!htb]
    \centering
    \begin{displaymath}
    \xymatrix{
        &    VGG^{n,1}_* \ar[r]\ar@/^1.8pc/@{.>}[rrr]^(.5){\scriptsize \txt{(vi)}}  &  WVGG_*^n \ar@{.>}[rr]^{\scriptsize \txt{(v)}} & & (SD^n)^C &\\
        &  & VGG^{n,n}_* \ar[u]   \\
        & VAG^n_* \ar[r]\ar[d] & VMG^n_* \ar[u]\ar[d]\\
        VG_*^n \ar[uuur] \ar@/_1.4pc/[rr] & WVAG^n_* \ar[r]  & WVMG^n_* \ar@/^-3pc/[uuu] \ar[]!<1.5pc,0pc>;[uuurr]
    }
    \end{displaymath}
    \caption{The non-self-decomposability of classes of strongly and weakly subordinated processes with arrows pointing in the direction of generalisation, assuming $n\!\ge\!2$,
        $|\Sigma|\neq 0$ and $\UUU((0,\infty)^n)>0$. Dotted arrows point in the direction of generalisation under additional assumptions corresponding to the parts of Theorem~\ref{thmVGGCnotseldec}.}
    \label{diagram3}
    %\end{center}
\end{figure}
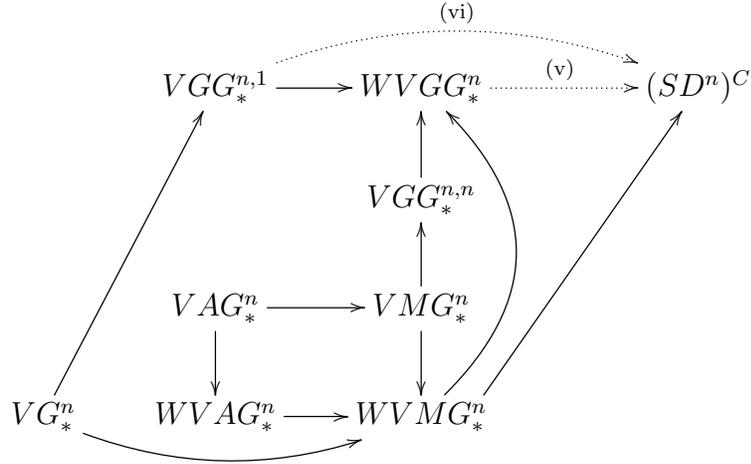

An $n$-dimensional subordinator $\bfT\!\sim\!AG^n_S(a,\bfalpha)$ is an {\em alpha-gamma subordinator}~\cite{Se08} provided
$\bfT\!=\!(T_1,\dots,T_n)\!\eqd\!G_0 \bfalpha\!+\!(G_1,\dots,G_n)$. A driftless Thorin subordinator $\bfT\!\sim\!\GGC^n_S({\bf 0},\UUU)$ is an alpha-gamma subordinator $\bfT\!\sim\!AG^n_S(a,\bfalpha)$ if and only if
$\UUU\!=\!a\bfdelta_{\bfalpha/\|\bfalpha\|^2}\!+\!\sum_{k=1}^n\beta_k\,\bfdelta_{\bfe_k/\alpha_k}$~(see~\cite[Remark~4.1]{BLMa}).\\[1mm]
\noindent{\bf Matrix gamma subordinator.}~Assume $\UUU$ is a finitely supported Thorin measure, that is~$\UUU\!=\!\sum_{k=1}^mu_k\,\bfdelta_{\bfalpha_k}$, where $u_k\!>\!0$, $\bfalpha_k\!\in\![0,\infty)^n_*$, $1\!\le\!k\!\le\!m$, $m\!\ge\!1$. An $n$-dimensional subordinator $\bfT\!\sim\!MG^n_S(\UUU)$ is a {\em matrix gamma subordinator} provided $\bfT\!\sim\!\GGC^n_S({\bf 0},\UUU)$~(see~\cite[Section~2.5]{BKMS16}).\\[1mm]
\noindent{\bf Strong variance alpha-gamma process.}~Assume $n\!\ge\!2$, $\bfB\!\sim\!BM^n(\bfmu,\allowbreak \Sigma)$ and $\bfT\!\sim\!AG^n_S(a,\bfalpha)$ are independent, where $\Sigma\!\in\![0,\infty)^{n\times n}$ a {\em diagonal} matrix. An $n$-dimensional process $\bfX\!\sim\!VAG^n(a,\bfalpha,\bfmu,\Sigma)$ is a {\em strong variance alpha-gamma process}~\cite{Se08} provided $\bfX\!\eqd\!\bfB\!\circ\!\bfT$.\\[1mm]
\noindent{\bf Weak variance alpha-gamma process.}~Assume $n\!\ge\!2$, $\bfB\!\sim\!BM^n(\bfmu,\allowbreak \Sigma)$ and $\bfT\!\sim\!AG^n_S(a,\bfalpha)$, where $\Sigma\!\in\!\RR^{n\times n}$ is an  {\em arbitrary} covariance matrix. An $n$-dimensional process $\bfX\!\sim\!WVAG^n(a,\bfalpha,\bfmu,\Sigma)$ is a {\em weak variance alpha-gamma process}~\cite{BLMa,BLMb,MiSz17} provided $\bfX\!\eqd\!\bfB\!\odot\!\bfT$.\\[1mm]
\noindent{\bf Strong variance matrix gamma process.}~ Assume $\bfB\!\sim\!BM^n(\bfmu,\allowbreak\Sigma)$ and $\bfT\!\sim\!MG^n_S(\UUU)$ are independent, where $\Sigma\!\in\![0,\infty)^{n\times n}$ a {\em diagonal} matrix. An $n$-dimensional process $\bfX\!\sim\!VMG^n(\bfmu,\Sigma,\UUU)$ is a {\em strong variance matrix gamma process} provided $\bfX\!\eqd\!\bfB\!\circ\!\bfT$.\\[1mm]
\noindent{\bf Weak variance matrix gamma process.}~Assume $\bfB\!\sim\!BM^n(\bfmu,\allowbreak\Sigma)$ and $\bfT\!\sim\!MG^n_S(\UUU)$, where $\Sigma\!\in\!\RR^{n\times n}$ is an  {\em arbitrary} covariance matrix. An $n$-dimensional process $\bfX\!\sim\!WVMG^n(\bfmu,\Sigma,\UUU)$ is a {\em weak variance matrix gamma process} provided $\bfX\!\eqd\!\bfB\!\odot\!\bfT$.\\[1mm]
\noindent{\bf Relation between subordination classes.}~The relationship between these classes of strongly and weakly subordinated processes is summarised in Figure~\ref{diagram1}.
Assume $n\!\geq\!2$. Theorem~\ref{thmVGGCseldec} and Figure~\ref{diagram1} give the chain of inclusions for self-decomposable processes depicted in Figure~\ref{diagram2}. Further assuming
$|\Sigma|\!\neq\!0$ and $\UUU((0,\infty)^n)\!>\!0$, Theorem~\ref{thmVGGCnotseldec} and Figure~\ref{diagram1} give the chain of inclusions for non-self-decomposable processes depicted in Figure~\ref{diagram3}.

Finally, using Theorem~\ref{thmVGGCseldec}(iii) and Theorem~\ref{thmVGGCnotseldec}(viii), we state the self-decomposability condition for $WVAG^n$ and $WVMG^n$ processes.
\begin{corollary}\label{corWValphaG}~Assume $\bfX\!\sim\!WVAG^n(a,\bfalpha,\bfmu,\Sigma)$.\\
    (i)~If $\bfmu\!=\!{\bf 0}$, then $\bfX\!\sim\!SD^n$.\\
    (ii)~If $|\Sigma|\!\neq\!0$ and $\bfmu\!\neq\!{\bf 0}$, then $\bfX\!\not\sim\!SD^n$.
\end{corollary}
\begin{corollary}\label{corWVMG}~Assume $\bfX\!\sim\!WVMG^n(\bfmu,\Sigma,\UUU)$.\\
    (i)~If $\bfmu\!=\!{\bf 0}$, then $\bfX\!\sim\!SD^n$.\\
    (ii)~If $|\Sigma|\!\neq\!0$, $\bfmu\!\neq\!{\bf 0}$ and $\UUU((0,\infty)^n)\!>\!0$, then $\bfX\!\not\sim\!SD^n$.
\end{corollary}
\begin{figure}[!htb]
    \centering
    \begin{displaymath}
    \xymatrix{
        VG^n_\bfnull\ar[rr]\ar[dddd]\ar[rd] &  &   WVGG^n_\bfnull \ar[d] \ar[r] & BM^n+GVGC_{\bfnull}^n\\
        &    VGG^{n,1} \ar[r]  &  WVGG^n & \\
        &  & VGG^{n,n} \ar[u]   \\
        & VAG^n \ar[r]\ar[d] & VMG^n \ar[u]\ar[d]\\
        VG^n \ar[uuur] \ar@/_1.4pc/[rr] & WVAG^n \ar[r]  & WVMG^n \ar@/^-3pc/[uuu]
    }
    \end{displaymath}
    \caption{The relationship between classes of strongly and weakly subordinated processes with arrows pointing in the direction of generalisation for $n\geq 1$ (note $VAG^n$ and $WVAG^n$ are not defined for $n=1$).}% The inclusion indicated by the arrow $\Rightarrow$ is proven in Theorem~\ref{thmVGGCseldec}, while the others are well-known (see \cite{BKMS16,BLMa}).
    \label{diagram1}
    %\end{center}
\end{figure}
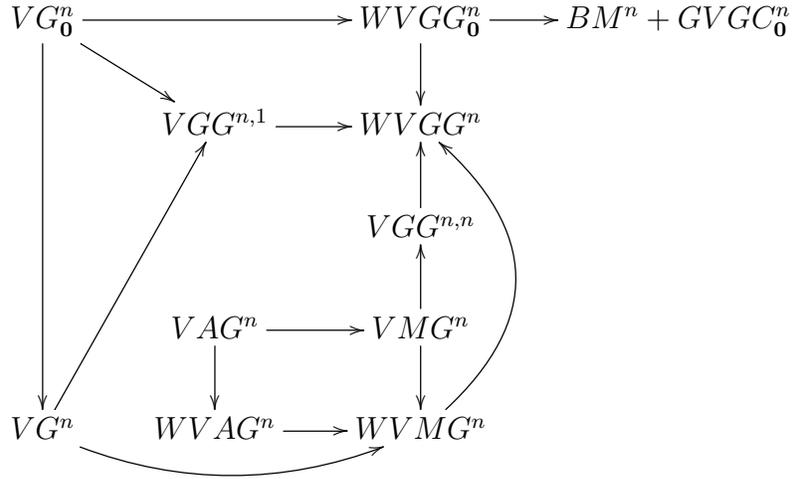
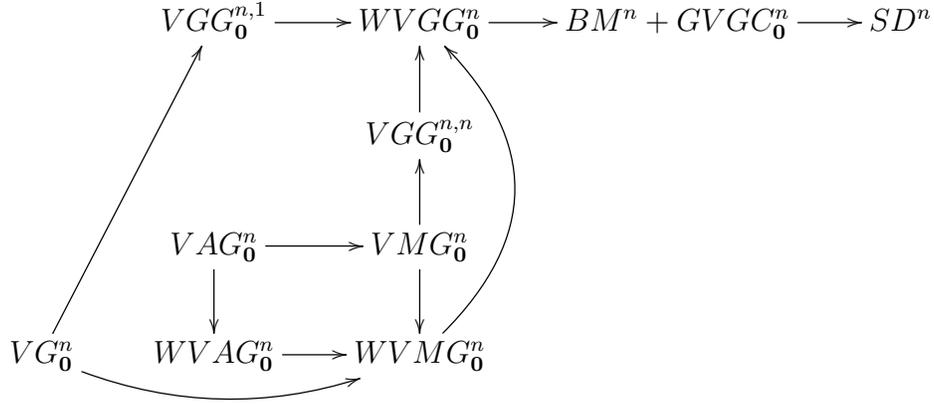
\begin{figure}[!htb]
    \centering
    \begin{displaymath}
    \xymatrix{
        &   VGG^{n,1}_\bfnull \ar[r]  &  WVGG_\bfnull^n \ar[r] & BM^n+GVGC_{\bfnull}^n\ar[r]&  SD^n\\
        &  & VGG^{n,n}_\bfnull \ar[u]   \\
        & VAG^n_\bfnull \ar[r]\ar[d] & VMG^n_\bfnull \ar[u]\ar[d]\\
        VG_\bfnull^n \ar[uuur] \ar@/_1.4pc/[rr] & WVAG^n_\bfnull \ar[r]  & WVMG^n_\bfnull \ar@/^-3pc/[uuu]
    }
    \end{displaymath}
\caption{The self-decomposability of classes of strongly and weakly subordinated processes with arrows pointing in the direction of generalisation, assuming $n\geq 2$.}
    \label{diagram2}
    %\end{center}
\end{figure}
\newpage\noindent{\bf Self-decomposable $\boldsymbol{WVGG_*}$ process.}~While Theorems~\ref{thmVGGCseldec} and \ref{thmVGGCnotseldec} seem to suggest that self-decomposability for $WVGG^n$ processes is equivalent to $\bfmu\!=\!\bfnull$, here we construct an example where this is not the case. For $n\!\geq\!2$, let $a,b,c\!\in\!(0,\infty)$, $g\!\in\![0,\infty)$, $\bfalpha\!\in\!(0,\infty)^n$, and $\UUU_k$, $1\!\le\!k\!\le\! n$, be univariate Thorin measures. Introduce the Thorin measure
\begin{equation}\label{sdcexthorin}
\UUU\;=\;\int_{g}^\infty \bfdelta_{u\bfalpha/\|\bfalpha\|^2}\frac{\rmd u}{(au+b)^{c}}\,+\,\sum_{k=1}^n\int_{(0,\infty)} \bfdelta_{u\bfe_k}\,\UUU_k(\rmd u)\,.
\end{equation}
Similar to an $AG^n$ subordinator, the subordinator $\bfT\!\sim\!GGC^n(\bfd,\UUU)$ associated to $\UUU$ has common and idiosyncratic time changes as it can be expressed as  $\bfT\!\eqd\!\bfd I\!+\!R_0\bfalpha\!+\!\sum_{k=1}^n R_k\bfe_k$, where $R_0\!\sim\!GGC^1_S(0,\eins_{(g,\infty)}(u)\rmd u/(au+b)^c)$, $R_k\!\sim\!GGC^1_S(0,\UUU_k)$, $1\!\le\!k\!\le\!n$.

Next, we provide an example of a $WVGG^n$ process with its Brownian motion subordinate having nonzero drift which is also self-decomposable, illustrating the sharpness of~\eqref{UUUAintegrab} (see Subsection~\ref{subsecproofpropsdcex} for a proof).
\begin{proposition}\label{propsdcex}
    Let $n\!\geq\!2$, $c\!\in\![1/2,1]$, $\bfd\!\in\![0,\infty)^n$, $\bfalpha\!\in\!(0,\infty)^n$, $\bfmu\!\in\!\RR^n_*$, $\Sigma\!\in\!\RR^{n\times n}$ be an invertible covariance matrix and $\UUU_k$, $1\!\le\! k\! \le\! n$, be univariate Thorin measures. Then there exist constants $a,b\!\in\!(0,\infty)$ and $g\!\in\![0,\infty)$ so that, with $\UUU$ defined in \eqref{sdcexthorin}, $\bfX\!\sim\!WVGG^n(\bfd,\bfmu,\Sigma,\UUU)$ is self-decomposable.
\end{proposition}
\noindent{\bf Thorin measures on the unit circle.}~Let $\bfX\!\sim\!WVGG^2(\bfd,\bfmu,\Sigma,\UUU)$ with $\bfmu\!\neq\!\bfnull$, $|\Sigma|\!\neq\!0$. If $\UUU_1\!=\!\int_{0}^1\bfdelta_{(\cos(\theta),\sin(\theta))}\,\rmd\theta$, then~\eqref{UUUAintegrabstronger} is satisfied. However, if $\UUU_2\!=\!\int_0^1\bfdelta_{(\cos(\theta^2),\sin(\theta^2))}\,\rmd\theta$,~\eqref{UUUAintegrabstronger} no longer holds as $\int_{0}^1 (\cos(\theta^2)\sin(\theta^2))^{-1/2}\,\rmd\theta$ is not finite. Both $\UUU_1$ and $\UUU_2$ are well-defined Thorin measure. Thus, $\bfX\!\not\sim\! SD^2$ if $\UUU\!=\!\UUU_1$, but no conclusion can be drawn on the basis of Theorem~\ref{thmVGGCnotseldec}(v) if $\UUU\!=\!\UUU_2$.

However, in the latter case, we show that \eqref{UUUAintegrab} holds. Let $\bfs\!=\!(s_1,s_2)\!\in\!\mySS_{**}$ and $\Sigma\!=\![\Sigma_{11},\Sigma_{12};\Sigma_{12},\Sigma_{22}]$. For $\theta\!\in\!(0,1]$, we have $\bfu\!=\!(\cos(\theta^2),\sin(\theta^2))\!\in\!(0,\infty)^2$ and $|\Sigma|\!\neq\!0$, which implies
\begin{align*}
\theta\mapsto f_\bfs(\theta):=\|\bfs\|^{-1}_{(\bfu\tr\Sigma)^{-1}}\big(\prod\bfu\big)^{-1/2}
\end{align*}
is continuous. When $\theta$ is small, we have
\begin{align*}f_\bfs(\theta)=\left(\frac{a_1\cos(\theta^2)+a_2\sin(\theta^2)}
{b_1\cos^2(\theta^2)+b_2\cos(\theta^2)\sin(\theta^2)} \right)^{1/2}
\end{align*}
for some constants $a_1,a_2,b_1,b_2\in\RR$. By noting $b_1\!=\!\Sigma_{11}s_2^2\!\neq\!0$ as $\bfs\!\in\!\mySS_{**}$, $|\Sigma|\!\neq\!0$, it follows that $\lim_{\theta\searrow0} f_\bfs(\theta)$ is finite, so the integral in \eqref{UUUAequivcirc}, which becomes $\int_0^1 f_\bfs(\theta)\,\rmd \theta$, is finite. Thus, Proposition~\ref{propUUUAequiv} implies that \eqref{UUUAintegrab} holds for all $\bfs\!\in\!\mySS_{**}$. Lastly, \eqref{meanpositvity} obviously holds on a subset of $\mySS_{**}$ having strictly positive Lebesgue surface measure by Theorem~\ref{thmpos}(v)--(vi) (see the proof of Theorem~\ref{thmVGGCnotseldec}(iv) for a similar argument). Thus, by Theorem~\ref{thmVGGCnotseldec}(iv), if $\UUU\!=\!\UUU_2$, then we also have $\bfX\!\not\sim\! SD^2$.\\[1mm]
\noindent{\bf $\boldsymbol{WVGG}$ processes from beta distributions of the second kind.}~A random variable $V(a,b)\!\eqd\!G_1/G_2$, where $G_1\!\sim\!\Gamma(a,1)$, $G_2\!\sim\!\Gamma(b,1)$ are independent, $a,b\!>\!0$, is said be \emph{beta distributed of the second kind.} Let $\UUU_{a,b}$ denote its probability measure. Here,
$\UUU_{a,b}(\rmd u)\!=\!C_{a,b}u^{a-1}(1\!+\!u)^{-a-b}\rmd u$, where $C_{a,b}$ is a normalising constant (see~\cite[Equation (2.2.5)]{Bo92}), is a Thorin measure satisfying \eqref{thorinmeasure} for all $a,b\!>\!0$.

Assume $n\!\geq\!2$. Let $\bfB\!\sim\!BM^n(\bfmu,\Sigma)$ and $\bfT\!=\!\sum_{k=1}^m T_k\bfalpha_k$, where $\bfalpha_k\!\in\![0,\infty)^n_*$, $1\!\le\! k\!\le\!m$, $m\!\geq\!1$, and $T_k\!\sim\!GGC_S^1(0,\UUU_{a_k,b_k})$, $a_k,b_k\!>\!0$, are independent for $1\!\le\! k\!\le\!m$. Let $\bfX \eqd \bfB\odot\bfT$.

By Theorem~\ref{thmVGGCseldec}, $\bfX\!\sim\! SD^n$ when $\bfmu\!=\!{\bf 0}$. With the additional assumption ${b_k}\!>\!1/2$, $1\!\le\! k \!\le\! m$, we have $\int_0^\infty u^{1/2}\,\UUU_{a_k,b_k}(\rmd u)\!<\!\infty$ since the integral is bounded near 0 by~\eqref{thorinmeasure} and bounded near infinity because, as $u\!\to\!\infty$, $u^{1/2}u^{a_k-1}(1\!+\!u)^{-a_k-b_k}\!\sim\!u^{-b_k-1/2}$. Thus, by Theorem~\ref{thmVGGCnotseldec}(vi), $\bfX \!\not\sim\! SD^n$ when $\bfmu\!\neq\!{\bf 0}$ and $|\Sigma|\!\neq\!0$. Note that this is an improvement to using the corresponding condition in Proposition~\ref{propGrig} as that would show non-self-decomposability only when $m\!=\!1$, $\bfalpha_1\!=\!\bfe$ with $b_1\!>\!1$.\\[1mm]
%\begin{remark}~\label{remCGMYGH}
\noindent{\bf $\boldsymbol{CGMY}$ and generalised hyperbolic processes.}~Multivariate $CGMY^n$ processes can be constructed as $VGG^{n,1}$ processes~\cite{MY08} or as $VGG^{n,n}$ processes~\cite{LS10}. Multivariate generalised hyperbolic processes $GH^n$ constructed as $VGG^{n,1}$ were considered in \cite[Example 1]{Gr07}. In all of these case, except for parameter choices where $GH^n$ reduces $VG^n$, if $\bfmu\!\neq \!\bfnull$, $|\Sigma|\!\neq\! 0$, then the integral in \eqref{defraymoments} is infinite, so we are unable to determine whether or not these processes are in $SD^n$. See~\cite[Examples 5.5.8 and 5.5.9]{Lu18} for additional details.

\begin{figure}[!htbp]
    \begin{center}
        \includegraphics{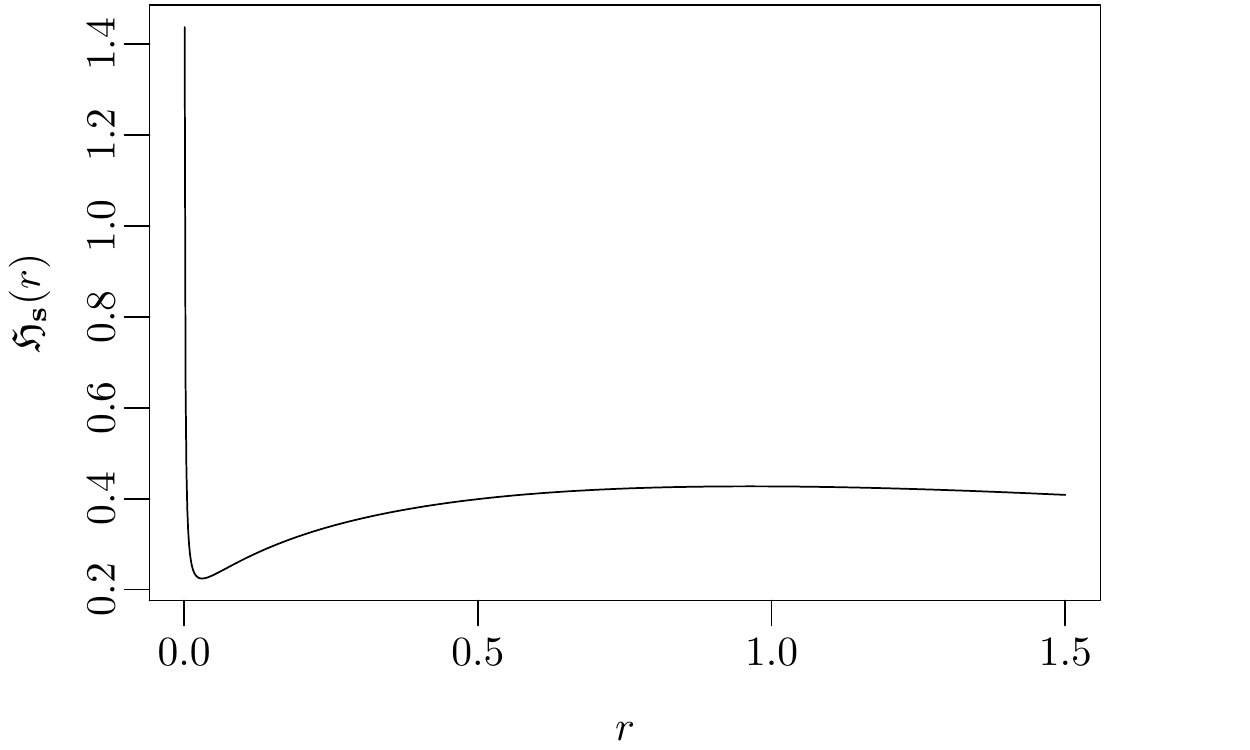}
        \caption{Plot of $r\mapsto \HHHH_{\bfs}(r)$ for $GH^3(-1,2,0.5,\bfmu,\diag(0.05,1,1))$, $\bfmu\!=\!(-5,1,1)$ at $\bfs\!=\!\bfmu/\|\bfmu\|$.}
        \label{plot1}
    \end{center}
\end{figure}

However, instead of using \eqref{defraymoments}, we can numerically examine the function $r\!\mapsto\!\HHHH_\bfs(r)$, $r\!>\!0$, $\bfs\!\in\!\mySS$, for particular parameter values. We use the parametrisation $GH^n(\alpha,\beta,\allowbreak\gamma,\bfmu,\Sigma)$ in \cite[Remark 2.3]{BKMS16}. Suppose that $\bfmu\!=\!(-5,1,1)$, $\bfX\!\sim\!GH^3(-1,2,0.5,\bfmu,\diag(0.05,\allowbreak1,1))$, a plot $r\!\mapsto\!\HHHH_\bfs(r)$ at $\bfs\!=\!\bfmu/\|\bfmu\|$ is given in Figure \ref{plot1}. If this behaviour extends to a set of $\bfs\!\in\!\mySS_{**}$ of strictly positive Lebesgue surface measure, $\bfX$ cannot be self-decomposable.\\[1mm]
%\halmos
%\end{remark}
\indent We conclude with the following remark.
\begin{remark}\label{remconclude}
    The proof of Theorem~\ref{thmVGGCnotseldec} attempts to show that $r\!\mapsto\!\HHHH_\bfs(r)$ is increasing at the origin to prove non-self-decomposability, while the numerical experiment above shows that it could be decreasing at the origin but strictly increasing at an alternative point, suggesting that to refine the conditions in Theorem~\ref{thmVGGCnotseldec} requires methods considering the entire domain.

    It would useful to determine sharper condition which could potentially be applied to multivariate $CGMY^n$ and $GH^n$ processes. Nevertheless, we believe that our analysis is an important step towards resolving this problem.

    Another possible research direction could be to extend our results about $WVGG^n$ processes to operator self-decomposability and to find conditions for their inclusion in Urbanik's $L$ classes. In the context of multivariate subordination, sufficient conditions for this were derived in \cite{BPS01}.
    \halmos
\end{remark}
\section{Proofs}\label{secproofs}
\subsection{Matrix Analysis for Proving Theorem~\ref{thmpos}(vi)}\label{subsectoolsunifpositivity}
We give four lemmas relating to the determinant, definiteness and invertibility of certain pattern matrices, which are needed to prove USP (Theorem~\ref{thmpos}(vi)).

To give an informal outline, let $\Sigma\!\in\!\RR^{n\times n}$ be an invertible covariance matrix. For some $\Delta\!\in\!\RR^{n\times n}$, we have  $\EEEE_{\bfmu,\Sigma}(\bfmu,\bfu)\!=\!\|\bfmu(\Delta\!*\!\Sigma)^{-1}\|_{(\Delta*\Sigma)^s}^2$, which is positive for $\bfu\!\in\!(0,\infty)^n$ provided $(\Delta\!*\!\Sigma)^s$ is positive definite, and moreover, this still holds under the limit approaching the infimum. To this end, Lemma~\ref{lemXi} shows that a relevant determinant is strictly positive, which implies in Lemma~\ref{lemUpsilon} that $\Upsilon\in\RR^{n\times n}$ is nonnegative definite. Consequently, Lemma~\ref{lemsymmutrSigma} establishes the positive definiteness of $2(\Delta*\Sigma)^s\!=\!\Upsilon*\Sigma$. Lemma~\ref{lemDelta} shows $\Delta\!*\!\Sigma$ is invertible, so that we can pass the result through the limit.

If $n\!=\!1$, set $\Xi_n(x)\!\equiv\!2$, $x\!\in\!\RR$, and otherwise, if $n\!=\!2,3,4\dots$
and $\bfx\!=\!(x_1,\dots,x_{n-1})\!\in\!\RR^{n-1}$, let $\Xi_n(\bfx)\!=\!(\Xi_{n,kl})\!\in\!\RR^{n\times n}$
be a square matrix, defined by
$\Xi_{n,kk}\!:=\!2$, $1\!\le\!k\!\le\!n$, $\Xi_{n,kl}\!=\!x_k$, $1\!\le\!k\!<\!n$, $l\!\neq\!k$, $\Xi_{n,nl}\!=\!1$, $1\!\le\!l\!<\!n$.
\begin{lemma}We have\label{lemXi} $\inf_{\bfx\in[0,1]^n}|\Xi_{n+1}(\bfx)|\!=\!2\!+\!n$ and
$\sup_{\bfx\in[0,1]^n}|\Xi_{n+1}(\bfx)|\!=\!2^{n+1}$ for $n\!\in\!\NN$.
\end{lemma}
\noindent{\bf Proof.}~If $n\!=\!1$, then this is clear. Otherwise, if $n\!\ge\!2$, note $\Xi_{n\!+\!1}(\bfe)$, $\bfe\!\in\!\RR^n$, turns into
the $(n\!+\!1)$-dimensional covariance matrix for equi-correlated variables, where $|\Xi_{n+1}(\bfe)|\!=\!2\!+\!n$~(see \cite[Theorem~8.4.4]{Gr83}).

If $n\!\ge\!1$, introduce a polynomial $h_n$ of degree $n$ in $\bfx\!=\!(x_1,\dots x_{n})\!\in\!\RR^{n}$
by $h_{n}(\bfx)\!:=\!|\Xi_{n\!+\!1}(\bfx)|$, $\bfx\!=\!(x_1,\dots x_{n})\!\in\!\RR^{n}$. Note $h_1(x_1)=4\!-\!x_1$, $h_2(x_1,x_2)=8\!-\!2x_1\!-\!2x_2$, $h_3(x_1,x_2,x_3)=16\!-\!4x_1-\!4x_2\!-\!4x_3+x_1x_2x_3$.

Expanding the determinants along its first row yields
$h_n(\bfx)\!=\!2h_{n-1}(\wt\bfx)\!+\!x_1r_{n-1}(\wt\bfx)$, where $\bfx\!=\!(x_1,\tilde\bfx)$, $x_1\!\in\!\RR$, $\wt\bfx\!=\!(x_2,\dots,x_{n})\!\in\!\RR^{n-1}$ and $\wt\bfx\!\mapsto\!r_{n-1}(\wt\bfx)$ is a remainder polynomial. Freezing
the variables $x_2,\dots,x_{n}$, $x\!\mapsto\!h_{n}(x,\wt\bfx)$ turns into an affine function so that $\partial_{x_1}^2h_{n}(x_1,\wt\bfx)\!\equiv\!0$. Note $h_n$ is invariant under coordinate permutations $h_n(\bfx P)\!=\!h_n(\bfx)$ for any per\-mutation matrix $P\in\RR^{n\times n}$. Particularly, $h_n$ is a harmonic function as div$(h_n)\!\equiv\!0$.

With $\partial[0,1]^n$ denoting the boundary of the set $[0,1]^n$ relative to $\RR^{n}$, the maximum principle for harmonic functions states that $\inf_{\bfx\in [0,1]^n}h_n(\bfx)\!=\!\min_{\bfx\in\partial[0,1]^n}h_n(\bfx)$
and $\sup_{\bfx\in [0,1]^n}h_n(\bfx)\!=\!\max_{\bfx\in \partial[0,1]^n}h_n(\bfx)$~(see \cite[Subsection~I.1.4]{Do01}).

Taking this into account and also the permutation invariance, observe that
 $\inf_{\bfx\in [0,1]^n}h_n(\bfx)\!=\!\min_{\bfx\in\myFF_n}h_n(\bfx)$
and $\sup_{\bfx\in [0,1]^n}h_n(\bfx)\!=\!\max_{\bfx\in\myFF_n}h_n(\bfx)$, where
\[\myFF_n\,:=\,\{{\bf0},\bfe\}\cup\bigcup_{1\le k < n}\{(x_1,\dots,x_n):
x_1=\dots=x_{k}=0,\,x_{k+1}=\dots =x_n=1\}\,.\]
For $n\!\ge\!2$, we have $h_n({\bf 0})\!=\!2^{n+1}$, $h_n({\bfe})\!=\!2\!+\!n$, while
\[h_n(\bfx)\,=\,2^k h_{n-k}(1,\dots, 1)\,=\,2^k(2\!+\!n\!-\!k),\;1\!\le\!k\!<\!n\,,\]
where $\bfx\!=\!(x_1,\dots,x_n)$, $x_1\!=\!\dots\!=\!x_{k}=0$, $x_{k+1}\!=\!\dots \!=\!x_n\!=\!1$, from which the result immediately follows.
\halmos\\[1mm]

\indent If $n\!\ge\!2$ and $\bfv\!=\!(v_1,\dots v_{n-1})\!\in\!\RR^{n-1}$, in\-tro\-duce a sym\-metric ma\-trix
$\Upsilon_n(\bfv)\!:=\!(\Upsilon_{n,kl}(\bfv))_{kl}\!\in\!\RR^{n\times n}$, where $\Upsilon_{n,kk}(\bfv)\!=\!2$, $1\!\le\!k\!\le\!n$ and $\Upsilon_{n,kl}(\bfv)\!=\!\Upsilon_{n,lk}(\bfv)\!=\!
1\!+\!\prod_{k\le m<l} v_{m}$, $1\!\le\!k\!<\!l\!\le\!n$.
\begin{lemma}\label{lemUpsilon} If $n\!\ge\!2$ and $\bfv\!\in\![0,1]^{n-1}$, then $\Upsilon_n(\bfv)$ is nonnegative definite.
\end{lemma}
\noindent{\bf Proof.}~Assume $n\!\ge\!2$ and $\bfv\!\in\![0,1]^{n-1}$.
Successively, from $k\!=\!1$ to $k\!=\!n\!-\!1$, in $\Upsilon_n(\bfv)$, multiply its $(k\!+\!1)$th column with $v_k$ and subtract
this from its $k$th column and then multiply its $(k\!+\!1)$th row with $v_k$ and subtract
this from its $k$th row. Afterwards, take out the factor
$x_k\!:=\!1\!-\!v_k$ from the $k$th column if $x_k\!\in\!(0,1]$, $1\!\le\!k\!<\!n$. If $x_k\!\in\!(0,1]$ for all $1\!\le\!k\!<\!n$, this yields $|\Upsilon_n(\bfv)|\!=\!|\Xi_n(\bfx)|\prod_{1\le k<n}x_k$,
where $\Xi_n(\bfx)=(\Xi_{n,kl}(\bfx))\in\RR^{n\times n}$ is the matrix in Lemma~\ref{lemXi} and thus, $|\Upsilon_n(\bfv)|\geq0$, $\bfv\in[0,1]^{n-1}$. Otherwise, if there is some $x_k=0$, $1\!\le\!k\!<\!n$, then the $k$th column is zero, so we get $|\Upsilon_n(\bfv)|=0$, $\bfv\!\in\![0,1]^{n-1}$.

Every other principal submatrix of $\Upsilon_n(\bfv)$, formed by keeping the rows and columns in the index set $\{a_1,\dots ,a_m\}$, $1\!\le\!a_1\!<\!\dots\!<\!a_m\!\le\!n$, $1\!\leq\!m\le\!n\!-\!1$, and the deleting the rest, is given by $\Upsilon_{m}(\overline{\bfv})$, where
\[\overline{\bfv}\,:=\,\Big(\prod_{k=a_1}^{a_2-1}v_k,\dots,\prod_{k=a_{m-1}}^{a_m-1}v_k\Big)\,\in\,[0,1]^{m-1}\,.\]
Using an argument similar to the preceding paragraph shows that $|\Upsilon_{m}(\overline{\bfv})|\!\ge\!0$. Therefore, $\Upsilon_n(\bfv)$ is nonnegative definite for all $\bfv\!\in\![0,1]^{n-1}$.~\halmos\\[1mm]
\indent If $\bfu\in(0,\infty)^n$ and $\Sigma\!\in\!\RR^{n\times n}$ is an invertible covariance matrix, introduce a symmetric matrix $\Sigma^s(\bfu,\Sigma)\in\RR^{n\times n}$ by
\begin{equation}\label{defSigmasymbfu}
    \Sigma^s(\bfu,\Sigma)\!:=\!((\bfu\tr\Sigma)\diag(1/\bfu))^s\,.
\end{equation}
\begin{lemma}\label{lemsymmutrSigma} Assume an invertible covariance matrix $\Sigma$. If $n\!\ge\!1$ and $\bfu\!\in\!(0,\infty)^n$, then $\Sigma^s(\bfu,\Sigma)$ is positive definite.\end{lemma}
\noindent{\bf Proof.} Introduce $t(u)\!:=\!(1\wedge u)\!+\!(1\wedge (1/u))\!\in\!(1,2]$, $u\!>\!0$, and the symmetric matrix $\Theta_n(\bfu)\!:=\!(\Theta_{n,kl}(\bfu))_{1\le k,l\le n}\!\in\!\RR^{n\times n}$ for $\Theta_{n,kl}(\bfu)\!:=\!t(u_k/u_l)$, $1\!\le\!k,l\!\le\!n,\bfu\!=\!(u_1,\dots,u_n)\!\in\!(0,\infty)^n$. By noting that $2\Sigma^s(\bfu,\Sigma)\!=\!\Theta_n(\bfu)\!*\!\Sigma$, $\bfu\!\in\!(0,\infty)^n$, $\Sigma$ is assumed to be an invertible covariance matrix, and $\Theta_n(\bfu)$ has positive diagonal entries $\Theta_{n,kk}(\bfu)\!\equiv\!2$, $1\!\le\!k\!\le\!n$, the positive definiteness of $\Sigma^s(\bfu,\Sigma)$ follows from
Oppenheim's inequality, provided we can show that
$\Theta_n(\bfu)$ is nonnegative definite, $n\!\in\!\NN$, $\bfu\!\in\!(0,\infty)^n$.

Indeed, $|\Theta_1(u)|\!\equiv\!2$, $u\!>\!0$, while $|\Theta_2(\bfu)|\!=\!4\!-\!t^2(u_1/u_2)\!\in\![0,3]$, $\bfu\!=\!(u_1,u_2)\!\in\!(0,\infty)^2$. Further, note $|\Theta_n(\bfu P)|\!=\!|P'\Theta_n(\bfu)P|\!=\!|\Theta_n(\bfu)|$, $\bfu\!\in\!(0,\infty)^n$, for per\-mu\-tation matrices $P\!\in\!\RR^{n\times n}$, $n\!\in\!\NN$.

Without loss of generality, we may assume
$\bfu\!\in\!(0,\infty)^n_{\le}$ and $n\!\ge\!2$. If $\bfu\!=\!(u_1,\dots,u_n)\!\in\!(0,\infty)^n_{\le}$, set $\bfv\!:=\!(v_1,\dots,v_{n-1})\!:=\!\bfv(\bfu)\!\in\!(0,1]^{n-1}$, for $v_k\!:=\!u_k/u_{k+1}$, $1\!\le\!k\!\le\!n$, and note $\Theta_n(\bfu)\!=\!\Upsilon_n(\bfv)$, where $\Upsilon_n(\bfv)$ is the nonnegative definite matrix in Lemma~\ref{lemUpsilon}, thus completing the proof.\halmos\\[1mm]
\indent If $n\!\ge\!2$ and $\bfv\!\in\!\RR^{n-1}$, introduce a matrix $\Delta_n(\bfv)\!:=\!(\Delta_{n,kl}(\bfv))\!\in\!\RR^{n\times n}$ and a subset $\GG_n\!\subseteq\!\RR^{n-1}$, where $\Delta_{n,kl}(\bfv)\!=\!1$, $1\!\le\!l\!\le\!k\!\le\!n$, $\Delta_{n,kl}(\bfv)\!=\!\Delta_{n,lk}(\bfv)\!=\!\prod_{k\le m<l} v_{m}$, $1\!\le\!k\!<\!l\!\le\!n$, and
\begin{equation}\label{gset}
\GG_n\,:=\,\bigcup_{k=1}^{n-1}\big\{\bfx=(x_1,\dots,x_{n-1})\in\partial[0,1]^{n-1}:x_k=0\big\}\,.
\end{equation}
Still assuming $n\!\ge\!2$, introduce a bijection from $(\mySS_{++})_\le$ to $(0,1]^{n-1}$ by
\begin{equation}\label{vtrafo}
\bfu\,=\,(u_1,\dots,u_n)\quad\mapsto\quad\bfv\,:=\,(u_1/u_2,\dots,u_{n-1}/u_n)\,.
\end{equation}
\begin{lemma}\label{lemDelta}Assume an invertible covariance matrix $\Sigma$. If $n\!\ge\!2$ and $\bfv\!\in\![0,1]^{n-1}$, then $\Delta_n(\bfv)*\Sigma$ is invertible.\end{lemma}
\noindent{\bf Proof.} Using the transformation $\bfu\!\mapsto\!\bfv=\bfv(\bfu)$ in~\eqref{vtrafo}, note $\Delta_n(\bfv)\!*\!\Sigma\!=\!(\bfu\!\tr\!\Sigma)\diag(1/\bfu)$ is invertible for $\bfv\!\in\!(0,1]^{n-1}$ because the RHS is, due to~\eqref{ineqOppHadutrSigma}. It remains to show invertibility for $\bfv\!\in\!\GG_n$.

If $n\!=\!2$, note $\GG_{n}\!=\!\{0\}$, so that $\Delta_2(\bfv)*\Sigma\!=\!\Delta_2(0)*\Sigma=[\Sigma_{11},0;\Sigma_{12},\Sigma_{22}]$ is indeed invertible. Next, assume $n\!\ge\!3$ and $\bfv\!=\!(v_1,\dots,v_{n-1})\!\in\!\GG_n$. If $v_{1}\!=\!0$ or $v_{n-1}\!=\!0$, note $\Delta_n(\bfv)*\Sigma=[\Sigma_{11},{\bf 0};\wt\bfsigma',\Delta_{n-1}(v_{2},\dots,v_{n-1})*\wt\Sigma]$
or $\Delta_n(\bfv)*\Sigma\!=\![\Delta_{n-1}(v_{1},\dots,v_{n-2})*\wt\Sigma,{\bf 0}';\wt\bfsigma,\Sigma_{nn}]$
for $\Sigma\!=\![\Sigma_{11},\wt\bfsigma;\wt\bfsigma',\wt\Sigma]$ and $\Sigma\!=\![\wt\Sigma,\wt\bfsigma';\wt\bfsigma,\Sigma_{nn}]$, respectively. Otherwise, we find a $1\!<\!k\!<n\!\!-\!1$ such that $v_k\!=\!0$, and
we may block $\Sigma$ into $\Sigma\!=\![\wt\Sigma_{11},\wt\Sigma_{12};\wt\Sigma_{21},\wt\Sigma_{22}]$
for $\wt\Sigma_{11}\!\in\!\RR^{k\times k},\wt\Sigma_{22}\!\in\!\RR^{(n-k)\times (n-k)}$ so that
\[\Delta_n(\bfv)*\Sigma\,=\,[\Delta_{k}(v_1,\dots,v_{k-1})*\wt\Sigma_{11},{0};\wt\Sigma_{21},
\Delta_{n-k}(v_{k+1},\dots,v_{n-1})*\wt\Sigma_{22}]\]
is invertible by induction hypothesis, completing the proof.
\halmos
\subsection{Proof of Theorem~\ref{thmpos}}\label{subsecproofofthmpos}
\noindent{\bf (i) and (ii).}~If $P\in\RR^{n\times n}$ is a permutation matrix, note
\begin{equation}\label{permutes}
\skal{\bfx}{\bfy}_{(\bfu \tr\Sigma)^{-1}}=\skal{\bfx P}{\bfy P}_{((\bfu P) \tr(P'\Sigma P))^{-1}}\,,\quad \bfx,\bfy\in\RR^n\,,\;\bfu\in(0,\infty)^n\,.
\end{equation}
If, in addition, $n\!\ge\!2$ and $\bfu\!\in\!(0,\infty)^n_{\le}$, we may write $\bfu\!=\!(\wt\bfu,u)$ for $u\!\in\!(0,\infty),\wt\bfu\!\in\!(0,\infty)^{n-1}_{\le}$ and $\Sigma\!=\![\wt\Sigma,\wt\bfsigma';\wt\bfsigma,\sigma]$ for $\sigma\!>\!0$, $\wt\bfsigma\!\in\!\RR^{n-1}$
and an invertible covariance matrix $\wt\Sigma\in\RR^{(n-1)\times(n-1)}$, thus having \\[-0.3cm]
\begin{equation}\label{redindim}\skal{\bfx}{\bfy}_{(\bfu \tr\Sigma)^{-1}}=\langle\wt\bfx,\wt\bfy\rangle_{(\wt\bfu\tr\wt\Sigma)^{-1}}+\frac{|\wt\bfu\tr\wt\Sigma|}{|\bfu\tr\Sigma|}\,
\big(x-\EEEE_{\wt\bfsigma,\wt\Sigma}(\wt\bfx,\wt\bfu)\big)\big(y-\EEEE_{\wt\bfsigma,\wt\Sigma}(\wt\bfy,\wt\bfu)
\big),\end{equation}
$\bfx\!=\!(\wt\bfx,x)$, $\bfy\!=\!(\wt\bfy,y)$, $\wt\bfx,\wt\bfy\!\in\! \RR^{n-1},x,y\!\in\!\RR$, with
$\EEEE_{\wt\bfsigma,\wt\Sigma}(\wt\bfx,\wt\bfu)\!=\!\langle\wt\bfx,\wt\bfu\tr\wt\bfsigma\rangle_{(\wt\bfu\tr\wt\Sigma)^{-1}}$ as defined in~\eqref{defADE}.

Recall $\bfy\!\in\!\RR_*^n$ and $\bfz\!\in\!\RR_{**}^n$, and set $i_\bfy\!:=\!\inf_{\bfu\in\mySS_{++}}\|\bfy\|^2_{(\bfu\tr\Sigma)^{-1}}$ and $I_\bfz\!:=\!\inf_{\bfu \in \mySS_{++}}\DDDD^2(\bfz,\bfu)$. Note $I_{\bfz}
\!=\!\inf_{\bfu\in(0,\infty)^n}\DDDD^2(\bfz,\bfu)$ using the invariance $\DDDD^2(\bfz,\bfu)\!=\!\DDDD^2(\bfz,\bfu^0)$ being true for $\bfu\!\in\!(0,\infty)^n_\le,\bfu^0\!:=\!\bfu/\|\bfu\|\!\in\!\mySS_{++}$.

If $n\!=\!1$ and $y,z\!\in\!\RR_*$, then the strict positivity of $i_y$ and $I_z$ is plain. Assume $n\!\ge\!2$. By choice of a suitable permutation in~\eqref{permutes} we may assume without loss of generality that
$\bfu\!=\!(\wt\bfu,u)\!=\!(u_1,\dots,u_n)\in(\mySS_{++})_\le,\wt\bfu\!\in\!(0,\infty)^{n-1}_{\le}$ and $u\!\in\!(0,1)$ with $u\!\ge\!\max\{u_k\!:\!1\!\le\!k\!<\!n\}$. Note $\|\wt\bfu\|\!\le\!\|\bfu\|\!=\!1$ and $\wt\bfu^0\!:=\!\wt\bfu/\|\wt\bfu\|\!\in\!(\mySS_{++})_\le$ while $u\!\geq\!n^{-1/2}$ because of $1\!=\!\sum_{k=1}^nu^2_k\!\le\!nu^2$.

If $\bfy\!=\!(\wt\bfy,y)\!\in\!\RR_*^n$, $\wt\bfy\!\in\!\RR^{n-1}$, $y\!\in\!\RR$,~\eqref{redindim} yields
$\|\bfy\|^2_{(\bfu \tr\Sigma)^{-1}}\!\ge\!\|\wt\bfy\|^2_{(\wt\bfu \tr\wt\Sigma)^{-1}}\!\ge\!\|\wt\bfy\|^2_{(\wt\bfu^{0}\tr\wt\Sigma)^{-1}}\!\ge\!i_{\wt\bfy}\!:=\!\inf_{\wt\bfu \in \mySS_{++}}\|\wt\bfy\|^2_{(\wt\bfu\tr\wt\Sigma)^{-1}}$; otherwise, if $\wt\bfy\!=\!{\bf 0}$ so that $y\!\neq\!0$,
combining~\eqref{ineqOppHadutrSigma} and~\eqref{redindim} yields
$\|\bfy\|^2_{(\bfu \tr\Sigma)^{-1}}\!\ge\!y^2|\wt\bfu\tr\wt\Sigma|/|\bfu\tr\Sigma|\!\ge\!y^2|\wt\Sigma|/\prod_{k=1}^n\Sigma_{kk}$ where the RHS does not depend on $\bfu$.

If $n\!\ge\!2$, suppose $\bfz\!=\!(\wt\bfz,z)\!\in\!\RR_{**}^n$, $\wt\bfz\!\in\!\RR_{**}^{n-1}$ and $z\!\in\!\RR_*$. Recalling
$\|\bfz\|^2_{(\bfu \tr\Sigma)^{-1}}\!\ge\!\|\wt\bfz\|^2_{(\wt\bfu \tr\wt\Sigma)^{-1}}$ and~\eqref{ineqOppHadutrSigma} yields, with $I_{\wt\bfz}\!:=\!\inf_{\widetilde\bfu \in \mySS_{++}}|\widetilde\bfu\tr\widetilde\Sigma|\|\widetilde\bfz\|^{2(n-1)}_{(\widetilde\bfu\tr\widetilde\Sigma)^{-1}}$,
\[
\DDDD^2_{\Sigma}(\bfz,\bfu)\ge\frac{|\bfu\tr\Sigma|}{|\wt\bfu^0\tr\wt\Sigma|}\;
\frac{\|\bfz\|^{2n}_{(\bfu\tr\Sigma)^{-1}}}{\|\wt\bfz\|^{2(n-1)}_{(\wt\bfu^0\tr\wt\Sigma)^{-1}}} I_{\wt\bfz}
=\frac{|\bfu\tr\Sigma|}{|\wt\bfu\tr\wt\Sigma|}\;
\frac{\|\bfz\|^{2n}_{(\bfu\tr\Sigma)^{-1}}}{\|\wt\bfz\|^{2(n-1)}_{(\wt\bfu\tr\wt\Sigma)^{-1}}} I_{\wt\bfz}\ge \frac{|\Sigma|}{n^{1/2}\prod_{k<n}\Sigma_{kk}}\,i_{\wt\bfz}\,I_{\wt\bfz}\,.\]
To summarise, strict positivity of $i_\bfy$ and $I_\bfz$ for all $n\!\in\!\NN$, $\bfy\!\in\!\RR^n_*$ and $\bfz\!\in\!\RR_{**}^n$ follows from mathematical induction.\\[1mm]
\noindent{\bf (iii) and (iv).}~Analogously, one shows that the suprema are finite~\cite[Lemma 5.4.4]{Lu18}.\\[1mm]
\noindent{\bf (v).}~Recall $\EEEE\!=\!\EEEE_{\bfmu,\Sigma}$
and $\VV^+\!=\!\VV_{\bfmu,\Sigma}^+$ in~\eqref{defADE} and~\eqref{defVV}, respectively.

It is straightforwardly verified that
$a\VV^+\!\subseteq\!\VV^+$, $a\!>\!0$,
and $\VV^+\!+\!\VV^+\!\subseteq\!\VV^+$
so that $\VV^+$ a convex cone.

The Cauchy-Schwarz inequality implies that
$|\EEEE(\bfx,\bfu)\!-\!\EEEE(\bfy,\bfu)|\!\le\!C\|\bfx\!-\!\bfy\|$ and thus
$\EEEE(\bfy,\bfu)\!\ge\!\EEEE(\bfx,\bfu)\!-\!C\|\bfx\!-\!\bfy\|$, $\bfx,\bfy\!\in\!\RR^n$, $\bfu\!\in\!(0,\infty)^n$, where $C\!:=\!1\!+\!\big(\sum_{k=1}^n\sup_{\bfu\in(0,\infty)^n}\EEEE^2(\bfe_k,\bfu)\big)^{1/2}$.
Note $C$ is a positive and finite constant, the latter by Part~(iv) of Theorem~\ref{thmpos}.
Assume $\bfx\!\in\!\VV^+$ so that $\underline \EEEE_{\bfx}\!:=\!\inf_{\bfu\in (0,\infty)^n}\EEEE(\bfx,\bfu)\!>\!0$. If $\bfy\!\in\!\RR^n$ so that $2C\|\bfx\!-\!\bfy\|\!\le\!\underline\EEEE_{\bfx}$, then noting $2\inf_{\bfu\in(0,\infty)^n}\EEEE(\bfy,\bfu)\!\ge\!\underline\EEEE_{\bfx}$ yields $\bfy\!\in\!\VV^+$, showing that $\VV^+$ is open.\\[1mm]
\noindent{\bf (vi).}~Assume $\bfmu\!\neq\!{\bf 0}$. We shall show that
$\bfmu\!\in\!\VV^+\!=\!\VV_{\bfmu,\Sigma}^+$.\\[1mm]
\noindent{\em Positivity.} Note~$\EEEE(\bfmu,\bfu)\!=\!\EEEE_{\mu,\Sigma}(\bfmu,\bfu)\!=\!\|\bfmu[(\bfu\tr\Sigma)\diag(1/\bfu)]^{-1}\|^2_{\Sigma^s(\bfu,\Sigma)}$, where $\Sigma^s(\bfu,\Sigma)\in\!\RR^{n\times n}$ in~\eqref{defSigmasymbfu} is symmetric and positive definite, the latter by Lemma~\ref{lemsymmutrSigma}, so that
$\EEEE(\bfmu,\bfu)\!>\!0$ for $\bfu\!\in\!(0,\infty)^n$ and $\bfmu\!\in\!\RR^n_*$.\\[1mm]
{\em Uniform strict positivity~(USP).}~If $n\!=\!1$, the result is obvious as $\EEEE(\mu,u)\!=\!\EEEE_{\mu,\Sigma}(\mu,u)\!\equiv\!\mu^2/\Sigma$, so that we may assume $n\!\ge\!2$ in the sequel. Let $\underline\EEEE_{\bfmu}\!:=\!\inf_{\bfu\in(0,\infty)^n}\EEEE(\bfmu,\bfu)$, and note
$\underline\EEEE_{\bfmu}\!=\!\inf_{\bfu\in\mySS_{++}}\EEEE(\bfmu,\bfu)$.
Plainly, there exists a sequence $(\bfu_m)_{m\ge 1}\!\subseteq\!\mySS_{++}$ such that $\lim_{m\to\infty}\EEEE(\bfmu,\bfu_m)\!=\!\underline\EEEE_{\bfmu}$.
Without loss of generality, choosing a suitable subsequence if necessary, we may assume that $\bfu_m\!\to\!\bfu_0$ for some $\bfu_0\!\in\!\mySS_+$, $m\!\to\!\infty$.

If $\bfu_0\!\in\!\mySS_{++}$, $\underline\EEEE_{\bfmu}\!=\!\EEEE(\bfmu,\bfu_0)\!>\!0$ follows from the
established positivity. Next, assume that $\bfu_0\!\in\!\mySS_+\backslash\mySS_{++}$ and $(\bfu_m)_{m\ge 1}\!\subseteq\!(0,\infty)^n_\le$, the latter without loss of generality by~\eqref{permutes}. Recall the definitions of $\GG_n$ and $\bfu\!\mapsto\!\bfv=\!\bfv(\bfu)$ in~\eqref{gset} and \eqref{vtrafo}, respectively. By selecting a suitable subsequence if necessary, we may assume that $\bfv_m\!:=\!\bfv(\bfu_m)\!\to\!\bfv_0\!\in\!\GG_n\!\subseteq\![0,1]^{n-1}$ as we assumed $\bfu_0\!\in\!\mySS_+\backslash\mySS_{++}$.

Note $\bfv\!\mapsto\!\Delta_n(\bfv)\!*\!\Sigma$ is continuous as a mapping from $\RR^{n-1}$ into $\RR^{n\times n}$. Since taking the inverse of matrices is a continuous operation, by applying Lemma~\ref{lemDelta} we must have $\EEEE(\bfmu,\bfu_m)\!=\!\bfmu(\Delta_n(\bfv_m)\!*\!\Sigma)^{-1}\bfmu'\!\to\!\bfmu(\Delta_n(\bfv_0)\!*\!\Sigma)^{-1}\bfmu'$ as $m\!\to\!\infty$, so that $\underline\EEEE_{\bfmu}\!=\!
\bfmu(\Delta_n(\bfv_0)\!*\!\Sigma)^{-1}\bfmu'$. As above, write
the RHS as $\underline\EEEE_{\bfmu}\!=\!\|\bfmu(\Delta_n(\bfv_0)\!*\!\Sigma)^{-1}\|^2_{(\Delta_n(\bfv_0)*\Sigma)^s}$.
Finally, note $2(\Delta_n(\bfv_0)\!*\!\Sigma)^s\!=\!\Upsilon_n(\bfv_0)\!*\!\Sigma$, where $\Upsilon_n(\bfv_0)$ is the nonnegative definite matrix in Lemma~\ref{lemUpsilon}, completing the proof of USP, that is~ $\underline\EEEE_{\bfmu}\!>\!0$. In particular,~$\bfmu\!\in\!\VV^+$ when $\bfmu\!\neq\!{\bf0 }$.~\halmos
\subsection{Proof of Theorem~\ref{thmHAEDformula}}\label{subsecproofthmHAEDformula}
Restricted to the nonnegative reals $r\!\in\!(0,\infty)$, $r\!\mapsto\!\KKKK_\rho(r)$ in~\eqref{besselsommerfeld} is nonincreasing for all $\rho\!\ge\!0$. If $\rho\!\ge\!0$, recall~(see~\cite[Equations~(8.451)--6 and (8.469)--3]{GrRy96})
\begin{equation}\label{Besselassinfty}
K_{\rho}(r)\,\sim\, K_{1/2}(r)\,=\,(\pi/2)^{1/2}\;e^{-r}/r^{1/2}\,,\quad r\to\infty\,.
\end{equation}
If $\rho\!>\!0$, $r\!\mapsto\!\KKKK_\rho(r)$ is uniformly bounded by $\KKKK_\rho(0+)\!=\!2^{\rho-1}\Gamma(\rho)$, as implication of~\eqref{besselsommerfeld}, while (see~\cite[Equation~(8.447)--3]{GrRy96})
\begin{equation}\label{Besselass0}
\KKKK_0(r)\,=\,K_{0}(r)\,\sim\,\ln^{-}r\,,\quad r \searrow 0\,.
\end{equation}
If $\rho\!\ge\!0$, we have
\begin{equation}\label{Besselunifbound}
\kappa_\rho\,:=\,\sup_{r>0}r\KKKK_\rho(r)\,<\,\infty\,,
\end{equation}
and
\begin{equation}\label{Bessellimo}
\lim_{r\searrow 0}r\KKKK_\rho(r)\,=\,0\,.
\end{equation}
If $\UUU$ is a Thorin measure on $[0,\infty)^n_*$, we get from~\eqref{thorinmeasure} that
\begin{equation}\label{hatBesselintegrabThorin}
\int_{(0,\infty)^n}\KKKK_{\rho}\big\{s\|\bfu\|^{\theta}\big\}\,\UUU(\rmd\bfu)\,<\,\infty\,,\quad s,\theta>0\,,\rho\ge 0\,.
\end{equation}
\begin{lemma}\label{lemprepHAEDformula} Assume a Thorin measure $\UUU$ on $(0,\infty)^n$.
If $f\!:\!(0,\infty)^n\!\to\!(0,\infty)$ is a Borel function such that, for some $a,b\!>\!0$, $f(\bfu)\!\ge\!a\|\bfu\|^{b}$
for all $\bfu\!\in\!(0,\infty)^n$, then
\begin{equation*}%\label{hatBesselintegrabThorin2}
\int_{(0,\infty)^n}\KKKK_{\rho}\{rf(\bfu)\}\,f^\theta(\bfu)\,\UUU(\rmd\bfu)\,<\,\infty\,, \qquad r>0,\;\rho,\theta\ge0\,.
\end{equation*}
\end{lemma}
\noindent{\bf Proof of Lemma~\ref{lemprepHAEDformula}.}~By~\eqref{Besselassinfty}, there exists $r_0\!=\!r_0(\rho,\theta)\!\in\!(0,\infty)$ such that
$\KKKK_{\rho}(r)\,r^\theta\!\le\!2\KKKK_{\rho+\theta}(r)$, $r\!>\!r_0$, implying
\[\int_{\{rf>r_0\}}\KKKK_{\rho}\{rf(\bfu)\}\,f^\theta(\bfu)\,
\UUU(\rmd\bfu)\,\le\, 2 r^{-\theta}\int_{(0,\infty)^n}\KKKK_{\rho+\theta}\{rf(\bfu)\}\,\UUU(\rmd\bfu)\,.\]
Since $\int_{(0,\infty)^n}\KKKK_{\rho+\theta}\{rf(\bfu)\}\,\UUU(\rmd\bfu)\!\le\!
\int_{(0,\infty)^n}\KKKK_{\rho+\theta}\{r a \|\bfu\|^b\}\,\UUU(\rmd\bfu)$,
the proof is completed
by~\eqref{hatBesselintegrabThorin} for large values of $\|\bfu\|$.

Set $h\!:=\!(r_0/ar)^{1/b}\!\in\!(0,\infty)$. If $\rho\!>\!0$, then
\[\int_{\{rf\le r_0\}}\KKKK_{\rho}\{rf(\bfu)\}\,f^\theta(\bfu)\,
\UUU(\rmd\bfu)\,\le\,(r_0/r)^\theta\KKKK_{\rho}(0+)\,\UUU\big(h\DD_{++}\big)\,.\]
If $\rho\!=\!0$,  recall $\KKKK_0\!\equiv\!K_0$, and note $c,d,g\!\in\![0,\infty)$, where $c\!:=\!\sup_{u>0}K_0(u)/(1\!+\!\ln^-u)$,
$d\!:=\!\sup_{u>0}(1\!+\!\ln^-\{rau^b\}))/(1\!+\!\ln^-u)$ and $g\!:=\!cd\,(r_0/r)^\theta$, so that
\[\int_{\{rf\le r_0\}}\KKKK_{0}\{rf(\bfu)\}\,f^\theta(\bfu)
\UUU(\rmd\bfu)\le g\int_{h\DD_{++}}
(1\!+\!\ln^-\|\bfu\|)\,\UUU(\rmd\bfu)\,<\,\infty\,.\]
Thus, the proof is completed for small values of $\|\bfu\|$.\halmos\\[1mm]
\noindent{\bf Proof of Theorem~\ref{thmHAEDformula}.}~Assume $\bfz\!\in\!\RR_{**}^n$ and $|\Sigma|\!\neq\!0$. Theorem~\ref{thmpos} states $\overline\EEEE_\bfz\!:=\!\sup_{\bfu\in(0,\infty)^n}|\EEEE(\bfz,\bfu)|\!<\!\infty$, $\zeta_\bfz\!:=\!\inf_{\bfu \in \mySS_{++}}\|\bfz\|_{(\bfu\tr\Sigma)^{-1}}\!>\!0$ and $\underline\DDDD_\bfz\!:=\!\inf_{\bfu\in(0,\infty)^n}\DDDD(\bfz,\bfu)\!>\!0$. The latter implies that $\UUU_\bfz(\rmd\bfu)\!:=\!\UUU(\rmd\bfu)/\DDDD(\bfz,\bfu)$ is a Thorin measure on $(0,\infty)^n$.\\[1mm]
{\bf (i).}~Setting $a\!:=\!2^{t/2}\zeta_\bfz^{t}\!>\!0$, we have $f(\bfu)\!:=\!\AAAA^t(\bfz,\bfu)\!\ge\!a\|\bfu\|^{t/2}$ so that Lemma~\ref{lemprepHAEDformula} states
\begin{equation}\label{hatBesselintegrabThorin3}
\int_{(0,\infty)^n}\KKKK_{\rho}\big\{r\AAAA^t(\bfz,\bfu)\big\}\AAAA^\theta(\bfz,\bfu)\,\UUU_\bfz(\rmd\bfu)\,<\,\infty\,,\qquad r,t>0,\;\rho,\theta\ge0.\end{equation}
If $w\!\in\!\mbox{dom}_\KKKK$, it is implied by~\eqref{besselsommerfeld} that $|\KKKK_{\rho}(w)|\!\le\!\KKKK_{\rho}\{\Re w^2\}$  so that, for $\bfu\!\in\!(0,\infty)^n$,
\[\big|\exp\{w\EEEE(\bfz,\bfu)\}\,\KKKK_{\rho}\big\{w\AAAA(\bfz,\bfu)\big\}\big|\,\le\,\exp\{\Re w\overline\EEEE_\bfz\}\,\KKKK_{\rho}\big\{\Re w^2\AAAA^2(\bfz,\bfu)\big\}\,.\]
In view of~\eqref{hatBesselintegrabThorin3}, the RHS and thus the LHS are $\UUU_\bfz(\rmd\bfu)$-integrable on $(0,\infty)^n$ in the last display. In particular,  $w\!\mapsto\!\HHHH_\bfz(w)$
in \eqref{defHHHH}, where $\rho\!=\!n/2$, is a well-defined function from $\mbox{dom}_\KKKK$ to $\CC$, while $\HHHH_\bfz((0,\infty))\!\subseteq\![0,\infty)$.

Next, let $\nu\!\ge\!1$ so that $\rho\!:=\!\nu\!-\!1\!\ge\!0$. We may differentiate~\eqref{besselsommerfeld} under the integral
to obtain $\rmd\KKKK_{\nu}(w)/\rmd w\!=\!-w\KKKK_{\rho}(w)$, $w\!\in\!\mbox{dom}_\KKKK$. In particular,
$w\!\mapsto\!e^{w\EEEE(\bfz,\bfu)}\,\KKKK_{\nu}\{w\AAAA(\bfz,\bfu)\}$ is holomorphic on $\mbox{dom}_\KKKK$ with derivative \[\EEEE(\bfz,\bfu)\exp\{w\EEEE(\bfz,\bfu)\}\,\KKKK_{\nu}\{w\AAAA(\bfz,\bfu)\}-w\AAAA^2(\bfz,\bfu)\exp\{w\EEEE(\bfz,\bfu)\}\,\KKKK_{\rho}\{w\AAAA(\bfz,\bfu)\}.\]
To see that $w\!\mapsto\!\HHHH_{\bfz}(w)$ can be differentiated under the integral, let $w\!\in\!\mbox{dom}_\KKKK$, and assume a compact and convex subset $K\!\subseteq\!\CC$ such that $\{0\}\!\subseteq\!K$ and $w\!+\!K\!\subseteq\!\mbox{dom}_\KKKK$. There exist $w_1,w_2,w_3\!\in\!w+K$ such that $\sup_{w'\in w+K}\Re w'\!=\!\Re w_1\!<\!\infty$, $\inf_{w'\in w+K}\Re(w')^2\!=\!\Re w^2_2\!>\!0$ and $\sup_{w'\in w+K} |w'|\!=\!|w_3|\!<\!\infty$,  so, for $\bfu\!\in\!(0,\infty)^n$,
\begin{equation*}\label{ineqproofHHH1}
\sup_{y\in K_*}\frac 1{|y|}\,\Big|\exp\{(w\!+\!y)\EEEE(\bfz,\bfu)\}\,\KKKK_{\nu}\{(w\!+\!y)\AAAA(\bfz,\bfu)\}-\exp\{w\EEEE(\bfz,\bfu)\}\,\KKKK_{\nu}\{w\AAAA(\bfz,\bfu)\}\Big|\end{equation*}
is bounded from above by
\begin{equation*}\label{ineqproofHHH2}\exp\{\overline\EEEE_\bfz\Re w_1\}\,\Big\{\overline\EEEE_\bfz \KKKK_{\nu}\{\Re w_2^2\AAAA^2(\bfz,\bfu)\}+|w_3|\,\AAAA^2(\bfz,\bfu)\,\KKKK_{\rho}\{\Re w_2^2\AAAA^2(\bfz,\bfu)\}\Big\}\,.\end{equation*}
Since $\rho\!=\!\nu\!-\!1\!\ge\!0$, it follows from~\eqref{hatBesselintegrabThorin3} that these expressions are $\UUU_\bfz(\rmd\bfu)$-integrable on $(0,\infty)^n$, completing the proof.\\[1mm]
\noindent{\bf (ii).}~If $\bfu\!\in\!(0,\infty)^n$, the Cauchy-Schwarz inequality states that $\EEEE^2(\bfz,\bfu)\!\le\!\AAAA^2(\bfz,\bfu)$.
In particular, the $\UUU_\bfz(\rmd\bfu)$-integrability of $\bfu\!\mapsto\!\AAAA(\bfz,\bfu)$ on $\bfu\!\in\!(0,\infty)^n$ is inherited by $\bfu\!\mapsto\!|\EEEE(\bfz,\bfu)|$.
Dominated convergence is applicable to both integrals in~\eqref{partialrHHHH}. The latter is verified by noting, for $\bfu\!\in\!(0,\infty)^n$,
\[\sup_{0<r\le 1}|\EEEE(\bfz,\bfu)\exp\{r\EEEE(\bfz,\bfu)\}
\KKKK_{n/2}\{r\AAAA(\bfz,\bfu)\}|\,\le\,\KKKK_{n/2}(0+)\exp\{\overline\EEEE_\bfz\}\,
\AAAA(\bfz,\bfu)\,,\]
and, by~\eqref{Besselunifbound} and $n\!\ge\!2$, that
\[\sup_{0<r\le 1}\exp\{r\EEEE(\bfz,\bfu)\}
r\AAAA^2(\bfz,\bfu)\,\KKKK_{(n-2)/2}\{r\AAAA(\bfz,\bfu)\}\,\le\,\kappa_{(n-2)/2}\,\exp\{\overline\EEEE_\bfz\}
\AAAA(\bfz,\bfu)\,.\]
This completes the proof of (ii) because the RHS in~\eqref{HHHHdiffat0} matches
\[\lim_{r\searrow 0}c_n\int_{(0,\infty)^n}\EEEE(\bfz,\bfu)\exp\{r\EEEE(\bfz,\bfu)\}
\KKKK_{n/2}\{r\AAAA(\bfz,\bfu)\}\,\UUU_\bfz(\rmd\bfu)\,,\]
while, by~\eqref{Bessellimo},
\[\lim_{r\searrow 0}c_n\int_{(0,\infty)^n}\exp\{r\EEEE(\bfz,\bfu)\}\,
r\AAAA^2(\bfz,\bfu)\,\KKKK_{(n-2)/2}\{r\AAAA(\bfz,\bfu)\}\,\UUU_\bfz(\rmd\bfu)=0\,.\]
\noindent{\bf (iii).}~Let $\bfX\!\sim\!W\VGGC^n(\bfd,\bfmu,\Sigma,\UUU)$. As shown by~\cite[Theorem~4.1]{BLMa}, the L\'evy measure $\XXX$ is related to the L\'evy measure $\VVV_{b,\bfmu,\Sigma}$ of the $VG^n(b,\bfmu,\Sigma)$ process as
\begin{equation*}
\XXX=\big\{(\UUU(\rmd \bfu)/\|\bfu\|^2)\otimes\VVV_{\|\bfu\|^2,\bfu\tr\bfmu,\bfu\tr\Sigma}(\rmd \bfx)\big\}\circ \big((\bfu,\bfx)\mapsto\bfx\big)^{-1}\,.
\end{equation*}
If, in addition, $|\Sigma|\!\neq\!0$, $\VVV_{b,\bfmu,\Sigma}$ admits a Lebesgue density $\nu_{b,\bfmu,\Sigma}\!=\!\rmd\VVV_{b,\bfmu,\Sigma}\big/\rmd\bfy$ (see~\cite[Equation~(2.11)]{BKMS16}), $\bfy\!\in\!\RR_*^n$,
\[\nu_{b,\bfmu,\Sigma}(\bfy)\,=\,
\frac{c_nb}{|\Sigma|^{1/2}\,\|\bfy\|^{n}_{\Sigma^{-1}}}\,
\exp\{\skal{\bfy}{\bfmu}_{\Sigma^{-1}}\}\,
\KKKK_{n/2}\big\{(2b\!+\!\|\bfmu\|^2_{\Sigma^{-1}})^{1/2}\|\bfy\|_{\Sigma^{-1}}\big\}\,.
\]
Let $\bfu\!\in\!(0,\infty)^n$, $\bfs\!\in\!\mySS_{**}$, $r\!>\!0$, $\bfz\!\in\!\RR_{**}^n$, and define
$\nu(\bfu,\bfz)\!:=\!\nu_{\|\bfu\|^2,\bfu\tr\bfmu,\bfu\tr\Sigma}(\bfz)$ so that
\[
\hhhh(\bfu,\bfs,r)\,:=\,c_n \exp\big\{r\EEEE(\bfs,\bfu)\big\}
\KKKK_{n/2}\big\{r \AAAA(\bfs,\bfu)\big\}\big/\DDDD(\bfs,\bfu)\,=\,r^{n}\nu(\bfu,r\bfs)\,.
\]
Combining these facts, the L\'evy measure $\XXX$ in Euclidean polar-coordinates of Borel sets $A\!\subseteq\!\RR_{**}^n$ is
\[
\XXX(A)\,=\,\int_{\mySS_{**}}\int_{0+}^{\infty}\int_{(0,\infty)^n} \eins_{A}(r\bfs)\,\hhhh(\bfu,\bfs,r) \,\UUU(\rmd \bfu)\,\frac{\rmd r}r\,\rmd\bfs\,.
\]
If $r\!\in\!(0,\infty)$, $\bfs\!\in\!\mySS_{**}$, integrating $\bfu\!\mapsto\!\hhhh(\bfu,\bfs,r)$ with respect to $\UUU(\rmd\bfu)$ on $(0,\infty)^n$ yields $\int_{(0,\infty)^n}\!\hhhh(\bfu,\bfs,r)\,\UUU(\rmd\bfu)\!=\!\HHHH_\bfs(r)$ in~\eqref{defHHHH}, as desired.~\halmos
\subsection{Proof of Theorem \ref{thmVGGCseldec}}
\label{subsecProofsthmselfdec}
If $n\!=\!1$, Theorem~\ref{thmVGGCseldec} specialises to Proposition~\ref{propGrig}(i). For the remaining part of the proof, suppose that $n\!\ge\!2$.\\[1mm]
{\bf (ii)}$\boldsymbol{\Rightarrow}${\bf (i).} Noting $VG^n_{\bf0}\!\subseteq\!\VGGC^{n,1}_{\bf0}$, Proposition~\ref{propGrig}(i) implies
that $VG^n_{{\bf0}}$ processes are $SD^n$ processes. Recall that the class of $SD^n$ processes is closed under convolution and convergence in distribution so that $GVGC_{\bf 0}\!\subseteq\!SD^n$. This completes the proof as adding in an independent Brownian motion does not effect the self-decomposability.\\[1mm]
{\bf (iii)}$\boldsymbol{\Rightarrow}${\bf (ii).} First assume $\bfX\!\sim\!W\VGGC^n_{{\bf 0}}({\bf0},\Sigma,\UUU)$ so that the characteristic exponent~\eqref{GVGcharexpo} simplifies to
\begin{equation}\label{dis0GVGcharexpo}
    \Psi_{\bfX}(\bftheta)\,=\,-\int_{[0,\infty)^n_*}\!\ln\big\{(\|\bfu\|^2
    +\frac 12\|\bftheta\|^2_{\bfu\tr\Sigma})\big/\|\bfu\|^2\big\}\,\UUU(\rmd\bfu)\,.
\end{equation}
We show that $\bfX\!\sim\!GVGC_{\bf0}^n$, the subclass of all distributions formed by convolutions and convergence in distribution of $VG^n_{\bf0}$ distributions.\\[1mm]
\noindent{\em Finitely supported Thorin measure.}~Suppose $\bfX\!\sim\!W\VGGC^{n}_{{\bf0}}({\bf 0},\Sigma,\UUU)$, where $\UUU\!=\!\sum_{k=1}^mu_k\bfdelta_{\bfalpha_k}$ for some $u_k\!\ge\!0$,
$\bfalpha_k\!\in\![0,\infty)^n_*$, and set $a_k\!:=\!\|\bfalpha_k\|$, $1\!\le\!k\!\le\!m$, $m\!\in\!\NN$ .

The characteristic exponent~\eqref{dis0GVGcharexpo} simplifies to
\[\Psi_\bfX(\bftheta)\,=\,-\sum_{k=1}^m u_k \ln\big\{(a_k^2\!+\!\frac 12\|\bftheta\|^2_{\bfalpha_k \tr\Sigma})/a_k^2\big\}\,,\quad\bftheta\!\in\!\RR^n\,.\]
Recall $\bfX\!\sim\!VG^n_{\bfnull}(b,\Sigma)$ has characteristic exponent (see~\cite[Equation~(2.9)]{BKMS16})
\begin{equation}\label{charVG}
\Psi_\bfX(\bftheta)\,=\,-b\ln\big\{(b+\frac 12\,\|\bftheta\|^2_\Sigma)/b\big\}\,,\quad \bftheta\in\RR^n\,.\end{equation}
Let $\bfX_1,\dots,\bfX_m$ be independent with
\[\bfX_k\sim VG^n_{\bf 0}(u_k,u_k\bfalpha_k\tr\Sigma/a_k^2)\,,\quad 1\!\le\! k\!\le\! m\,.\]
Then $\Psi_\bfX\!=\!\sum_{k=1}^m\Psi_{\bfX_k}$ as~\eqref{charVG} implies
\begin{equation*}\Psi_{\bfX_k}(\bftheta)%&=&-u_k\ln\{(\|w_k\bfe/\|\bfe\|^2\|^2+\frac 12\|\bftheta\|^2_{(w_k\bfe/\|\bfe\|^2)\tr(\bfalpha_k\tr\Sigma/w_k)})/\|w_k\bfe/\|\bfe\|^2\|^2\}\\
    =-u_k\ln\{(a_k^2+\frac 12\|\bftheta\|^2_{\bfalpha_k\tr\Sigma})/a_k^2\}\,,\quad \bftheta\!\in\!\RR^n\,,\quad 1\!\le\! k\!\le\! m\,.
\end{equation*}
{\em Arbitrary nonzero Thorin measure.}~Let $\bfX\!\sim\!WVGG^{n}_{{\bf 0}}({\bf 0},\Sigma,\UUU)$, where $\UUU$ is an arbitrary nonzero Thorin measure.

Introduce $w(r)\!:=\!(1\!+\!\ln^-r)\!\wedge\!(1/r)$, $r\!>\!0$ and, for $\bftheta\!\in\!\RR^n$, a nonnegative and continuous function
\begin{equation}\bfu\mapsto g_{\bftheta}(\bfu)\,:=\,\ln\{(\|\bfu\|^2+\frac{1}{2}\|\bftheta\|^2_{\bfu\tr\Sigma})/\|\bfu\|^2\}\big/w(\|\bfu\|),\quad \bfu\in[0,\infty)^n_*\,.\label{defgtheta}\end{equation}
As $\mySS_+$ is compact and $\bfu\!\mapsto\!\|\bftheta\|^2_{\bfu\tr\Sigma}$ is continuous on $[0,\infty)^n_*$, $2c\!:=\!\sup_{\bfs\in \mySS_+}\|\bftheta\|^2_{\bfs\tr\Sigma}$ is a finite constant.
If $\bfu\!\in\![0,\infty)^n_*$, note $0\!\le\!g_{\bftheta}(\bfu)\!\le\!\widetilde{g}(\|\bfu\|)$, with $\widetilde{g}(r)\!:=\!\ln\{1\!+\!(c/r)\}\big/w(r)$, $r\!>\!0$. Here $\widetilde{g}(r)$ can be continuously extended by setting $\widetilde{g}(0)\!:=\!\widetilde{g}(0+)\!=\!1$. Besides this, note $\widetilde{g}(r)\!=\!r\ln\{1+(c/r)\}\!\le\!c$, $r\!\ge\!1$. Thus, $\widetilde{g}$ is uniformly bounded on $[0,\infty)$ and so is $g_{\bftheta}$ on $[0,\infty)^n_*$. To summarise, $g_{\bftheta}$ is continuous and bounded on its domain $[0,\infty)^n_*$ for all $\bftheta\!\in\!\RR^n$.

As a Thorin measure, $\UUU$ satisfies~\eqref{thorinmeasure} so that~$I\!:=\!\int_{[0,\infty)^n_*}w(\|\bfu\|)\,\UUU(\rmd\bfu)\!\in\!(0,\infty)$ as we assumed $\UUU$ is nonzero measure. We introduce a Borel probability measure on $[0,\infty)^n_*$ by $\PPP(\rmd\bfu)\!:=\!w(\|\bfu\|)\UUU(\rmd\bfu)/I$.

Perceiving ${\bf 0}$ and ${\bf \infty}$ as infinitely far points,
$[0,\infty)^n_*$ is a locally compact space, so there exists a sequence of
finitely supported Borel probability measures $(\PPP_k)_{k\in\NN}$ such that $\PPP_k$ converges weakly to $\PPP$ as $k\!\to\!\infty$ (see
\cite[Corollaries 30.5 and 30.9]{Ba01}). For $k\!\in\!\NN$, introduce a finitely supported Thorin measure $\UUU_k(\rmd\bfu):=I\PPP_k(\rmd\bfu)/w(\|\bfu\|)$ with an associated process $\bfX_k\!\sim\! W\VGGC^n_{{\bf0}}({\bf0},\Sigma,\UUU_k)$.

Let $\bftheta\!\in\!\RR^n$, and recall $g_{\bftheta}$ in~\eqref{defgtheta} is continuous and bounded so that
\[
\Psi_{\bfX_k}(\bftheta)\!=\!
-I\int_{[0,\infty)_*^n} g_{\bftheta}(\bfu)\PPP_k(\rmd \bfu)%\\
\!\to\!-I\int_{[0,\infty)_*^n} g_{\bftheta}(\bfu)\PPP(\rmd \bfu)
\!=\!\Psi_{\bfX}(\bftheta)\,,\quad k\!\to\!\infty\,.\]
This proves that for any $\bfX\!\sim\!WVGG^{n}_{{\bf 0}}({\bf 0},\Sigma,\UUU)$ with nonzero $\UUU$, there exists a sequence of processes $(\bfX_k)_{k\in\NN}$, where $\bfX_k$, $k\!\in\!\NN$, is the
superposition of finitely many independent $VG^n_{\bf0}$ processes, and $\bfX_k\!\to\!\bfX$ in law as $k\to\infty$.\\[1mm]
{\em General case}. Assume $\bfX\!\sim\!WVGG^{n}_{{\bf0}}(\bfd,\Sigma,\UUU)$ with no restrictions on $\bfd$, $\Sigma$ and $\UUU$. If $\UUU\!=\!0$, then $\bfX\!\sim\!WVGG^{n}_{{\bf 0}}({\bfd},\Sigma,\UUU)\!=\!BM^n({\bf 0},\bfd\!\tr\!\Sigma)$ by~\eqref{GVGcharexpo}, which is in the required form as the null process ${\bf0}\!\sim\! VG^{n}_{{\bf0}}(b,0)$ is an element of $GVGC^n_{\bf0}$. The case where $\UUU\!\neq\!0$ and $\bfd\!=\!{\bf 0}$ is dealt with above.
If $\UUU\!\neq\!0$ and $\bfd\!\neq\!{\bf 0}$, then $\bfX\!\eqd\!\bfB\!+\!\bfY$ by~\eqref{GVGcharexpo} where $\bfB\sim BM^n({\bf 0},\bfd\tr\Sigma)$ and $\bfY
\!\sim\!WVGG^{n}_{{\bf 0}}({\bf 0},\Sigma,\UUU)$ are independent processes. Since we showed above that $\bfY\sim GVGC^n_{\bf0}$, this completes the proof. \halmos
\begin{remark}
    Assume $\bfd\!=\!\bfnull$. In the proof of Theorem~\ref{thmVGGCseldec}, we show that $W\VGGC_\bfnull^n\!\subseteq\!GVGC_\bfnull^n$. It is natural to ask whether $W\VGGC_\bfnull^n\!=\!GVGC_\bfnull^n$. Clearly, this holds if and only if the $W\VGGC^n_\bfnull$ class is closed under convolution and convergence in distribution, which is true when $n\!=\!1$ (see \cite[Theorem~7.3.1]{Bo92}). But for $n\!\geq\!2$, we conjecture $W\VGGC^n_\bfnull$ distributions are not closed under convolutions. The $n\!=\!1$ result relies on the fact that the parameter $\Sigma$ can be absorbed into the subordinator, whereas this cannot easily be done when $n\!\geq\!2$, which explains the difficulty in extending the $n\!=\!1$ result and why we expect $W\VGGC_\bfnull^n\!=\!GVGC_\bfnull^n$ to be false in the latter case.\halmos
\end{remark}
\subsection{Proof of Theorem~\ref{thmVGGCnotseldec}}\label{subsecProofthmsnotseldec}
{\bf (ii)}$\boldsymbol{\Rightarrow}${\bf (i).}~This is a modification of the arguments in~\cite[Pro\-po\-sitions~1 and~3(ii)]{Gr07}.\\[1mm]
{\bf (iii)}$\boldsymbol{\Rightarrow}${\bf (ii).}~This is obvious.\\[1mm]
{\bf (iv)}$\boldsymbol{\Rightarrow}${\bf (iii).}~This follows from Theorem~\ref{thmHAEDformula}(ii), see~\eqref{HHHHdiffat0}.\\[1mm]
{\bf (v)}$\boldsymbol{\Rightarrow}${\bf (iv).}~Theorem~\ref{thmpos}(v) states that
$\VV^+\!=\!\VV^+_{\bfmu,\Sigma}$ in~\eqref{defVV} is an open convex cone. As we assumed $\bfmu\!\neq\!{\bf 0}$, note
$\bfmu\!\in\!\VV^+$ by Theorem~\ref{thmpos}(vi) so that $\VV^+$ is not the empty set. In particular, $\VV^+_{**}$ is a nonempty open convex cone
so that $\BB^+\!:=\!\mySS\!\cap\!\VV^+_{**}$ is a nonempty and relative-open subset of $\mySS_{**}$, thus having
strictly positive $(n\!-\!1)$-dimensional Le\-besgue surface mea\-sure.

Let $\bfs\!\in\!\BB^+$, and note $\underline\EEEE_\bfs\!=\!\inf_{\bfu\in(0,\infty)^n}\EEEE(\bfs,\bfu)\!>\!0$ as $\bfs\!\in\!\VV^+$ so that
\[\int_{(0,\infty)^n}\,\EEEE(\bfs,\bfu)\,\frac{\UUU(\rmd\bfu)}{\DDDD(\bfs,\bfu)}\;\ge\;\underline\EEEE_\bfs \int_{(0,\infty)^n}\,\frac{\UUU(\rmd\bfu)}{\DDDD(\bfs,\bfu)}\,.\]
Also, noting $\{\bfu\!\in\!(0,\infty)^n\!:\DDDD(\bfs,\bfu)>0\}\!=\!(0,\infty)^n$, the RHS in the last display is strictly positive as we assumed $\UUU((0,\infty)^n)\!>\!0$, thus~\eqref{meanpositvity} holds.

Theorem~\ref{thmpos}(iii) states the finiteness of $\xi\!:=\!\sup_{\bfu\in\mySS_{++}}\|\bfu\tr\bfmu\|^2_{(\bfu\tr\Sigma)^{-1}}$ while
Theorem~\ref{thmpos}(i) states strict positivity of $\zeta_\bfs\!:=\!\inf_{\bfu\in \mySS_{++}}\|\bfs\|^2_{(\bfu\tr \Sigma)^{-1}}$ since $\bfs\in\BB^+\!\subseteq\!\RR_*^n$. If $\bfu\!\in\!(0,\infty)^n$ and $\bfu^0\!:=\!\bfu/\|\bfu\|$, noting
$\|\bfu\|^{n}|\bfu^0\tr\Sigma|\!=\!|\bfu\tr\Sigma|$ and using~\eqref{ineqOppHadutrSigma} yields
\[\frac{\AAAA^2(\bfs,\bfu)}{\DDDD^2(\bfs,\bfu)}\,=\,\frac{2\|\bfu\|\!+\!\|\bfu^0\tr\bfmu\|^2_{(\bfu^0\tr\Sigma)^{-1}}}
{\|\bfs\|^{2(n-1)}_{(\bfu^0\tr\Sigma)^{-1}}|\bfu^0\tr\Sigma|}\,\le\,
\frac{2\vee \xi}{|\Sigma|\zeta_\bfs^{n\!-\!1}}
\frac{(1\!+\!\|\bfu\|)\|\bfu\|^n}{\prod\bfu}\,.\]
Particularly,~\eqref{UUUAintegrab} holds for $\bfs\!\in\!\BB^+$ as we assumed~\eqref{UUUAintegrabstronger}.\\[1mm]
{\bf (vi)}$\boldsymbol{\Rightarrow}${\bf (v).}~If $1\!\le\! k\!\le\! m$, set $u^*_k\!:=\!\prod\bfu_k$, and note
%\begin{align}\label{raysubordsimpf}
%\|v\bfu_k\|^n\!=\!v^n\|\bfu_k\|^n\,,\quad \prod(v\bfu_k)\!=\!v^nu_k^*\,,\quad v>0,\;\bfu_k\in(0,\infty)^n,
%\end{align}
$\|v\bfu_k\|^n\!=\!v^n\|\bfu_k\|^n$ and $\prod(v\bfu_k)\!=\!v^nu_k^*$, $v\!>\!0$,
so that
\[\int_{(0,\infty)^n}(1\!+\!\|\bfu\|^{1/2})\frac{\|\bfu\|^{n/2}\,\UUU(\rmd\bfu)}{(\prod\bfu)^{1/2}}=
\sum_{k=1}^m\left(\frac{\|\bfu_k\|^n}{u^*_k}\right)^{1/2}\int_{(0,\infty)}(1+\|v\bfu_k\|^{1/2})\,\UUU_k(\rmd v)\,.\]
Note the RHS in the last display is bounded above by
\[\sum_{k=1}^m\left(\frac{\|\bfu_k\|^n}{u^*_k}\right)^{1/2}(1\vee\|\bfu_k\|^{1/2})\,\int_{0+}^\infty (1+v^{1/2})\,\UUU_k(\rmd v)\,<\,\infty\]
since $\UUU_k((0,1])<\infty$, $1\!\le\!k\!\le\!n$. Thus,~\eqref{UUUAintegrabstronger} holds, completing the proof of the implication.\\[1mm]
{\bf (vii)}$\boldsymbol{\Rightarrow}${\bf (vi).}~Plainly, finitely supported Thorin measures satisfy~\eqref{defraymoments}.~\halmos
\subsection{Proof of Proposition~\ref{propUUUAequiv}}\label{subsecProofUAequi}
Let $\bfmu\!\in\!\RR^n$. By recalling $\xi\!:=\!\sup_{\bfu\in\mySS_{++}}\|\bfu\tr\bfmu\|^2_{(\bfu\tr\Sigma)^{-1}}$ is finite due to Theorem~\ref{thmpos}(iii), and using \eqref{ineqOppHadutrSigma}, we have (see the proof of Theorem~\ref{thmVGGCnotseldec}(v) for a similar argument)
\begin{equation}\label{adequivboundpf}
\frac{2}{\prod_{k=1}^n\Sigma_{kk}}\frac{\|\bfu\|^2}{\|\bfs\|^{2(n-1)}_{(\bfu\tr\Sigma)^{-1}}\prod\bfu} \le \frac{\AAAA^2(\bfs,\bfu)}{\DDDD^2(\bfs,\bfu)} \le \frac{2\vee \xi}{|\Sigma|}\frac{(1+\|\bfu\|)\|\bfu\|}{\|\bfs\|^{2(n-1)}_{(\bfu\tr\Sigma)^{-1}}\prod\bfu}\,.
\end{equation}
{\bf (i).}~By noting $(1\!+\!r)r\!\leq\!(1\!+\!1/a)r^2$ for $r\!>\!a$, $a\!>\!  0$, it follows from \eqref{adequivboundpf} that \eqref{UUUAintegrab} is equivalent to \eqref{UUUAequiv}.\\[1mm]
{\bf (ii).}~This is immediately follows from~\eqref{adequivboundpf} by noting $\|\bfu\|\!=\!1$ for $\bfu\!\in\!\mySS_{++}$.\\[1mm]
{\bf (iii).}~By substituting $\bfu\!=\!v\bfu_{k}$, $v\!>\!0$, $\bfu_k\!\in\!(0,\infty)^n$, $1\!\le\!k\!\le\!m$, $m\!\le\!1$, in~\eqref{adequivboundpf}, it follows from~\eqref{thorinmeasure} that~\eqref{UUUAintegrab} is equivalent to \eqref{defraymomentsequiv} (see the proof of Theorem~\ref{thmVGGCnotseldec}(vi) for a similar argument). \halmos
\begin{remark}\label{remarkwhynogenadequiv}
Note $\inf_{\bfu\in\mySS_{++}}\|\bfu\tr\bfmu\|^2_{(\bfu\tr\Sigma)^{-1}}$ can be 0 (see Remark~\ref{remconverseinfsup}), so the term $\|\bfu\|^2$ in the LHS of \eqref{adequivboundpf} cannot be replaced by a constant multiple of $(1+\|\bfu\|)\|\bfu\|$. Thus, we are unable to obtain an equivalent condition to \eqref{UUUAintegrab} for general Thorin measures.\halmos
\end{remark}
\subsection{Proof of Proposition~\ref{propsdcex}}\label{subsecproofpropsdcex}
Let $n$, $c$, $\bfd$, $\bfalpha$, $\bfmu$, $\Sigma$ and $\UUU_k$, $1\!\le\!k\!\le\! n$, be as specified in the assertion of the Proposition~\ref{propsdcex}. Set $\nu\!:=\!n/2$ so that $\rho\!:=\!\nu\!-\!1\ge 0$.

Note $\bfalpha\!\tr\!\bfmu\neq\bfnull$ and, thus, $\|\bfalpha\tr\bfmu\|_{(\bfalpha\tr\Sigma)^{-1}}\!>\!0$ so that we find $a,b\!\in\!(0,\infty)$, satisfying
\begin{equation}\label{countereq}
2b\,=\,a\|\bfalpha\tr\bfmu\|^2_{(\bfalpha\tr\Sigma)^{-1}}\,.\end{equation}
Let $\bfs\!\in\!\mySS$. By compactness of $\mySS$ and continuity of
$\bfs\!\mapsto\!a_\bfs\!:=\!2\|\bfs\|^2_{(\bfalpha\tr\Sigma)^{-1}}$, $\bfs\!\mapsto\! b_\bfs\!:=\!\|\bfalpha\!\tr\!\bfmu\|^2_{(\bfalpha\tr\Sigma)^{-1}}\|\bfs\|^2_{(\bfalpha\tr\Sigma)^{-1}}$ and $\bfs\!\mapsto\!\EEEE_\bfs\!:=\!\skal{\bfs}{\bfalpha\tr\bfmu}_{(\bfalpha\tr\Sigma)^{-1}}$, we have $\underline{a}\!:=\!\inf_{\bfs\in\mySS}a_\bfs\!>\!0$, $\underline{b}\!:=\!\inf_{\bfs\in\mySS}b_\bfs\!>\!0$ and $\overline{\EEEE}\!:=\!\sup_{\bfs\in\mySS}\EEEE_{\bfs}\!\in\!(0,\infty)$, respectively. Finally, set $h\!:=\!\overline{\EEEE}\pi^{1/2}\Gamma(\nu\!+\!1/2)/\Gamma(\nu)$ and $g\!:=\!1\vee((h^2\!-\!\underline{b})/\underline{a})$.

It is straightforwardly checked that $\UUU^*\!:=\!\int_{g}^\infty \bfdelta_{u\bfalpha/\|\bfalpha\|^2}\rmd u/(au+b)^{c}$ is a Thorin measure, and so is $\UUU$ in~\eqref{sdcexthorin}, the latter being a finite superposition of Thorin measures. Besides this, there are independent processes $\bfB$, $\bfX^*$, $X_1,\dots X_n$, where $\bfB\!\sim\!BM^n(\bfd\tr\bfmu,\bfd\tr\Sigma)$, $X_k\!\sim\!\VGGC^{1,1}(0,\mu_k,\Sigma_{kk},\UUU_k)$, $1\!\le\!k\!\le\!n$, and $\bfX^*\!\sim\!WVGG^{n}(\bfnull,\bfmu,\Sigma,\UUU^*)$ so that $\bfX\!\eqd\!\bfB\!+\!\bfX^*\!+\!\sum_{k=1}^nX_k\bfe_k$ (see~\cite[Proposition~3.7 and Remark~3.9]{BLMa}). Recall $\bfB\!\sim\!SD^n$, $X_k\!\sim\!SD^1$, $1\!\le\!k\!\le\!n$, by Theorem~\ref{thmVGGCseldec}.

Self-decomposability is closed under convolution, so it suffices to show that $\bfX^*\!\sim\!SD^n$. The latter holds, provided we can show that, for all $\bfs\!\in\!\mySS$, $r\!\mapsto\!\HHHH^*_{\bfs}(r)$ in nondecreasing, where $\HHHH^*_{\bfs}(r)$ is~\eqref{defHHHH} applied to $\bfX^*$. For $r\!>\!0$, $\bfs\!\in\!\mySS$, we have
\[
\HHHH^*_{\bfs}(r)\,:=\, \frac{2}{(2\pi)^{\nu}}\,\frac{\exp\{r\EEEE_\bfs\}}{\|\bfs\|^n_{(\bfalpha\tr\Sigma)^{-1}}|\bfalpha\tr\Sigma|^{1/2}}\,\int_{g}^\infty\KKKK_{\nu}\{r(a_\bfs u\!+\!b_\bfs)^{1/2}\}\frac{\rmd u}{({a}u\!+\!{b})^c}\,.
\]
Recall~\eqref{countereq}, and note ${b}-{a}b_\bfs/a_\bfs\!\equiv\!0$. Also, note~$\rmd\KKKK_{\nu}(r)/\rmd r\!=\!-r\KKKK_{\rho}(r)$. We can differentiate $r\!\mapsto\!\HHHH^*_{\bfs}(r)$ using the dominated convergence theorem and Lemma~\ref{lemprepHAEDformula} (see the proof of Theorem~\ref{thmHAEDformula}(i) for a similar argument), then making the change of variables $v\!=\!r(a_\bfs u\!+\!b_\bfs)^{1/2}$, we obtain for $\bfs\!\in\!\mySS$, $r\!>\!0$,
\begin{eqnarray}\label{defHHHHstar}
\frac{\partial_r \HHHH^*_{\bfs}(r)}{\HHHH^*_{\bfs}(r)}&=&\EEEE_\bfs-r\frac{\int_{g}^{\infty}(a_\bfs u\!+\!b_\bfs)\KKKK_{\rho}\{r(a_\bfs u\!+\!b_\bfs)^{1/2}\}\rmd u/({a}u\!+\!{b})^c}{\int_{g}^{\infty}\KKKK_{\nu}\{r(a_\bfs u\!+ \!b_\bfs)^{1/2}\}\rmd u/({a}u\!+\!{b})^c}\nonumber\\
&=&\EEEE_\bfs- \frac{\int_{r(a_\bfs g\!+\!b_\bfs)^{1/2}}^{\infty}\KKKK_{\rho}(v)v^{3-2c} \,\rmd v}{r\int_{r(a_\bfs g+b_\bfs)^{1/2}}^{\infty}\KKKK_{\nu}(v)v^{1-2c}\,\rmd v}\,.\label{nonSDcexcond}
\end{eqnarray}
Note \cite[Equation~(2.11)]{Gau14} states that
\begin{align*}
\int_{r}^\infty\KKKK_{\nu}(v)\,\rmd v\,\le\,\frac{\pi^{1/2}\Gamma(\nu\!+\!1/2)}{\Gamma(\nu)}\,\KKKK_{\nu}(r),\quad r>0\,,
\end{align*}
and, as $1\!-\!2c\!\le\!0$, this inequality yields
\[
\overline{\EEEE}\int_{r(a_\bfs g+b_\bfs)^{1/2}}^{\infty}\KKKK_{\nu}(v)v^{1-2c}\,\rmd v\,\le\, h r^{1-2c}(a_\bfs g+b_\bfs)^{(1-2c)/2}\,\KKKK_{\nu}\{r(a_\bfs g+b_\bfs)^{1/2}\}\,.\label{nonSDcexi1}
\]
Since $2\!-\!2c\!\ge\!0$, we have
\[
\int_{r(a_\bfs g+b_\bfs)^{1/2}}^{\infty}\KKKK_{\rho}(v)v^{3-2c}\,\ge\, r^{2-2c}(a_\bfs g+b_\bfs)^{1-c}\int_{r(a_\bfs g+b_\bfs)^{1/2}}^{\infty} v\KKKK_{\rho}(v)\,\rmd u\,,
\]
and, by the fundamental theorem of calculus with $\rmd\KKKK_{\nu}(v)/\rmd v\!=\!-v\KKKK_{\rho}(v)$,
\[
\frac{1}{r}\int_{r(a_\bfs g+b_\bfs)^{1/2}}^{\infty}\KKKK_{\rho}(v)v^{3-2c}\,\rmd v\,\ge\,r^{1-2c}(a_\bfs g+b_\bfs)^{1-c}\KKKK_{\nu}\{r(a_\bfs g+b_\bfs)^{1/2}\}\,.\label{nonSDcexi2}
\]
Recall $\underline{a}g\!+\!\underline{b}\!\ge\!h^2$ by choice of $g$, so that combining the above results, we have
$h\partial_r\ln\HHHH^*_{\bfs}(r)\!\le\!\overline{\EEEE}(h\!-\!(\underline{a}g\!+\!\underline{b})^{1/2})\!\le\!0$, $\bfs\!\in\!\mySS$,
completing the proof.\halmos
\begin{remark}
    In the situation of Proposition~\ref{propsdcex}, using the notation in the proof above, we show that the truncation parameter $g$ can be reduced to $0$ when $n\!=\!2$ and $c\!=\!1/2$.
    Introduce $t\!:=\!rb_{\bfs}^{1/2}\!>\!0$. Let $\bfs\!\in\!\mySS$, $f_\bfs\!:=\!\EEEE_{\bfs}/b_{\bfs}^{1/2}$, $g\!=\!0$. It follows from integration by parts that
    $\partial_r\ln\HHHH^*_\bfs(r)\!\le\!0$ in~\eqref{defHHHHstar} is equivalent to $\mathfrak{G}^*_\bfs(t)\!\ge\!0$, where
    \begin{align*}
    \mathfrak{G}^*_\bfs(t):=\KKKK_{1}(t) +(1/t-f_\bfs)\int_{t}^\infty \KKKK_{1}(v)\,\rmd v.
    \end{align*}
    For $t\!>\!0$, taking the derivative yields
    \begin{align*}\partial_t\mathfrak{G}^*_\bfs(t)=&-t\KKKK_{0}(t)+(f_\bfs - 1/t)\KKKK_1 (t)-t^{-2}\int_{t}^\infty \KKKK_{1}(v)\,\rmd v\\
    &\leq(f_\bfs t - 1)K_1(t)-tK_0(t)\,.
    \end{align*}
    Thus, $\partial_t\mathfrak{G}^*_\bfs(t)\!<\!0$ for all $t\!>\!0$ since
    \begin{equation*}%\label{besselratiobound}
    f_\bfs-\frac{1}{t}\,\leq\,1-\frac{1}{t}\,<\,\frac{t}{1+(1+t^2)^{1/2}}\,<\,\frac{K_{0}(t)}{K_1(t)},\quad t>0\,,
    \end{equation*}
    where the first inequality follows from the Cauchy-Schwarz inequality, the second inequality is obvious and the third inequality follows from~\cite[Theorem~1.2]{LN10}.
    %    \begin{align*}
    %    1-\frac{d}{t} <\frac{t}{\nu+(\nu^2+t^2)^{1/2}},\quad t>0\,,
    %    \end{align*}
    So $t\!\mapsto \!\mathfrak{G}^*_\bfs(t)$, $t\!>\!0$, has no turning points. In addition, we have that $\lim_{t\searrow 0}\mathfrak{G}^*_\bfs(t)\!=\!\infty$ and $\lim_{t\to \infty}\mathfrak{G}^*_\bfs(t)\!=\!0$. Thus, $\bfX\!\sim\! SD^2$ for $g\!=\!0$ .

    If $n\!=\!2$ and $c\!\in\!(1/2,1]$, or $n\!\geq \!3$ and $c\!\in\![1/2,1]$, we leave as as an open problem whether the truncation parameter can be chosen to be $g\!=\!0$ while $\bfX\!\not\sim \!SD^n$.
    %still producing examples of self-decomposable processes with Brownian nonzero drift

    In \eqref{UUUAintegrab}, the integral is a multiple of $\int_g^\infty (a_\bfs u\!+\!b_\bfs)^{1/2}\,\rmd u/(au+b)^c\!=\!\infty$ for all $\bfs\!\in\!\mySS$ and $c\!\in\![1/2,1]$, so Proposition~\ref{propsdcex} is not in contradiction with Theorem~\ref{thmVGGCnotseldec}.\halmos
\end{remark}

\section*{Acknowledgement}
We thank the editor and the two referee for their suggestions that have substantially improved the paper. K.W.~Lu's research was supported by an Australian Government Research Training Program Scholarship. B.~Buchmann and K.W.~Lu's research was supported by ARC grant DP160104737.

\end{document}